\documentclass[letterpaper, 10 pt, conference]{ieeeconf}
\IEEEoverridecommandlockouts  
\usepackage{amsmath} 
\usepackage{amssymb}  
\usepackage{amsthm,graphicx,epstopdf}
\usepackage[usenames, dvipsnames]{xcolor}
\usepackage{soul}
\usepackage{subcaption}
\usepackage{mwe}
\usepackage{array}
\usepackage{balance}
\usepackage{esvect}
\usepackage{url}
\usepackage{float}
\usepackage{bm}
\usepackage{pgf,tikz}\usepackage{mathrsfs}\usetikzlibrary{arrows}
\usepackage{multirow}
\usepackage{bigstrut}
\usepackage{empheq}
\usepackage{hyperref}
\usepackage{enumerate}
\usepackage{romannum}
\usepackage{tcolorbox}

\definecolor{LightGray}{rgb}{.95,.95,.95}

\theoremstyle{plain}
\newtheorem{thm}{Theorem}

\theoremstyle{remark}
\newtheorem{remark}{Remark:}
\theoremstyle{remark}

\theoremstyle{remark}

\theoremstyle{remark}

\DeclareMathAlphabet{\pazocal}{OMS}{zplm}{m}{n}

\graphicspath{{./figures/}}

\newlength\dlf  

\newcommand*\widefbox[1]{\fbox{\hspace{2em}#1\hspace{2em}}}

\newcommand{\hatmap}[1]{{#1}^{\times}}
\newcommand{\hatmapname}{cross-map}
\newcommand{\T}[1]{\tilde{#1}}
\newcommand{\Om}{\Omega}

\newcommand{\dvone}{\tilde{\eta}}
\newcommand{\dvtwo}{\tilde{\theta}}

\newcommand{\norm}[1]{\left\lVert#1\right\rVert}
\newcommand{\Rho}{\mathrm{P}}

\newcommand{\Jtrue}{\bar{J}}
\newcommand{\Jnom}{J}

\newcommand{\qmodel}{\textit{true model}}
\newcommand{\Qmodel}{\textit{True model}}

\newcommand{\refsys}{\textit{reference model}}
\newcommand{\Refsys}{\textit{Reference model}}

\newcommand{\highlight}[1]{{#1}} 
\newcommand{\rebuttal}[1]{#1}

\title{\Large \textbf{Geometric $L_1$ Adaptive Attitude Control for a Quadrotor Unmanned Aerial Vehicle}}

\author{Prasanth Kotaru$^{1}$, Ryan Edmonson$^{2}$, and Koushil Sreenath$^{1}$\\
	\thanks{This work is supported in part by NSF Grant CMMI-1840219, PITA, Autel, and in part by the Google faculty research award.}
	\thanks{$^{1}$P. Kotaru and K. Sreenath are with the Dept. of Mechanical Engineering, University of California, Berkeley, CA, 94720 email: {\tt \{ prasanth.kotaru, koushils\}@berkeley.edu}
		.}
	\thanks{$^{2}$R. Edmonson is with Calspan, Buffalo, NY, email:
		{\tt ryan.edmonson@calspan.com}.}
}

\begin{document}
	\setlength{\abovedisplayskip}{2pt}
	\setlength{\belowdisplayskip}{2pt}
	\setlength{\textfloatsep}{5pt}	
	
	\maketitle
	
	\begin{abstract}
In this paper, we study the quadrotor  UAV attitude
control on $SO(3)$ in the presence of unknown disturbances
and model uncertainties. $L_1$ adaptive control for UAVs using Euler angles/quaternions is shown to exhibit robustness and precise attitude tracking in the presence of disturbances and uncertainties. \highlight{However, it is well known that dynamical models and controllers that use Euler angle representations are prone to singularities and typically have smaller regions of attraction while quaternion representations are subject to the unwinding phenomenon.}
To avoid such complexities, we present a Geometric $L_1$ adaptation
control law to estimate the uncertainties. A model reference
adaptive control approach is implemented, with the attitude errors between the quadrotor model and the reference model defined on the manifold. 
Control laws for the quadrotor and reference models are developed directly on $SO(3)$ to track the desired trajectory while rejecting the uncertainties. 
Control Lyapunov function based analysis is used to show the exponential input-to-state
stability of the attitude errors. The proposed $L_1$ adaptive
controller is validated using numerical simulations. Preliminary
experimental results are shown comparing a geometric PD
controller to the geometric $L_1$ adaptive controller. Experimental validation of the proposed controller is carried out on an Autel X-star quadrotor. 
\end{abstract}
	
	\section{Introduction}
\label{sec: intro}

In recent years, quadrotor unmanned aerial vehicles (UAVs) have been an area of increasing interest. Due to their small size and simple mechanical structure, quadrotors have a large range of potential applications including visual inspection and transportation, as well as a medium for testing control techniques for research purposes.
Generally, an attitude controller for a quadrotor uses Euler angles or quaternions as the attitude states of the system. Instead of these typical attitude controllers, in geometric control, the entries of the rotation matrix between the body-fixed frame and inertial frame are used as the attitude states. Geometric control can be used for complex flight maneuvers as seen in \cite{lee2011geo} and it completely avoids singularities and complexities that arise when using local coordinates. Geometric control has also been used in robust tracking \cite{lee2011geo} and carrying suspended loads with cables \cite{koushil2013geo, lee2015adapt, wu2014geo}.

\begin{figure}
	\centering
	\includegraphics[width=0.8\columnwidth]{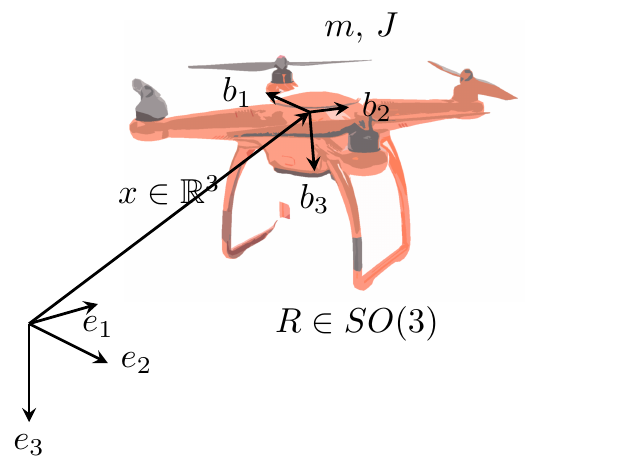}
	\caption{An Autel X-star quadrotor system is shown. The Quadrotor model evolves on $SE(3)$, with the attitude represented by a rotation matrix $R\in SO(3)$, from body-frame to inertial-frame. The origin of the body-frame $\bm{\{b_1, b_2, b_3\}}$ is at the center-of-mass of the quadrotor. \highlight{Video link for the experiments: \href{https://youtu.be/nBDDxpkz6Pg}{https://youtu.be/nBDDxpkz6Pg} }} 
	\label{fig:quad_geo}
\end{figure}

However, these controllers are model dependent and assume accurate representation of the dynamics. In the presence of disturbances or uncertainties, the controllers would result in unstable control or large tracking errors. Lately, there has been much work done with adaptive control, to achieve higher performance and robustness in the presence of model uncertainties and disturbances. 

\subsection{\rebuttal{Related Work}}
Model Reference Adaptive Control (MRAC) is a standard approach to adaptive control, where an ideal reference system is simulated, and the differences between the real
and reference systems are used to predict and cancel the
disturbances in the system \cite{whitehead2010mrac, ghaffar2016mrac}. 
However, this approach has some important practical drawbacks. If the adaptation rate is pushed too high, the system is given high-frequency input commands, which may not be feasible for the system. \highlight{Also, high adaptation gain might result in deteriorated closed-loop stability or unpredictable transient behaviours \cite{zuo2014adapt, cao2006design}.} Furthermore, the high-frequency commands may excite the system, which is not desired. Because of this, the adaptation rate will be low, leading to a long prediction in the system  uncertainties. 
$L_{1}$ adaptive control is an extension of MRAC, where the adaptation and control are decoupled through the use of a low pass filter. The designer can choose the cutoff frequency for the low pass filter, removing the high-frequency components of the adaptation from the control input. This promotes a very fast adaptation rate while keeping the control input sufficiently smooth, making $L_{1}$ { a much more practical adaptive control technique \cite{dodenhoft2017design, ackerman2016l1}}.

$L_1$ attitude controllers have been developed for attitude control of quadrotors, using Euler angles (or quaternions) as attitude states  \cite{zuo2014adapt, michini2009adapt, monte2013adapt, mallikarjunan2012l1} or for general linear systems \cite{Cao2006adapt}. Geometric adaptive schemes have been previously developed for quadrotor control in \cite{lee2015adapt, kulumani2016adapt}. Similar to \cite{lee2015adapt}, our work also develops the attitude dynamics and control laws for the quadrotor directly on the $SO(3)$ manifold without any local simplifications. However, in contrast to \cite{lee2015adapt}, we implement a model reference approach, consider time-varying disturbances and uncertainties, estimate uncertainties using a $L_1$ adaptation scheme, and show exponential input-to-state stability. 

\subsection{\rebuttal{Challenges}}
\rebuttal{Challenges when implementing the $L_1$ adaptation scheme to a geometric setting is to ensure correct formulation for the errors between the \refsys~(or state predictor) and \qmodel~(or plant). A simple difference between the states of the state predictor and the plant cannot be used\cite{hovakimyan2010l1}, because these states belong to different tangent spaces and thus their difference often lead to incorrect measure of the error, especially for large errors. Similarly, care should be taken when defining the uncertainty prediction error. The adaptation law for the uncertainty requires a careful choice of prediction law that ensures the Lyapunov function for the prediction errors decreases with time.}

\subsection{\rebuttal{Contribution}}
In this paper, we study the problem of a quadrotor with unknown disturbances and model uncertainties. 

\highlight{The main focus of this work is on implementing the $L_1$ adaptation directly on the $SO(3)$ manifold. In particular, the contributions of this paper with respect to prior work are,}
\begin{enumerate}
	\item We define attitude tracking errors, between the quadrotor model and a reference model (without disturbances), directly on the tangent bundle of the $SO(3)$ manifold and develop the control law to show that these attitude errors exhibit exponential input-to-state stability.
	\item We propose a geometric $L_1$ adaptation law to estimate the uncertainties and formally prove the resulting geometric input-to-state stability.
	\item We present numerical results to validate the performance of the proposed geometric $L_1$ adaptive control and show that the proposed control outperforms \highlight{geometric control \textit{without} $L_1$ (geometric $PD$ control in \cite{lee2011geo})}. 
	
	\item \highlight{We also compare the geometric $L_1$ with a traditional Euler $L_1$ and numerically establish that the Geometric $L_1$ outperforms the Euler $L_1$.}
    \item Finally, we present experimental results to show the tracking performance of the $L_1$ adaptive control w.r.t. the baseline geometric PD controller (\textit{without} $L_1$).
\end{enumerate}

\subsection{\rebuttal{Organization}}
The rest of the paper is structured as follows. Section~\ref{section: quad_dynamics} revisits the dynamics and control of quadrotor on SE(3). Section~\ref{sec:disturb} defines dynamics in the presence of disturbances and model uncertainties and presents the dynamics in terms of attitude tracking errors. Section~\ref{section: control} proposes a reference model without any uncertainties and defines attitude errors between the actual model and reference model. It proposes the Geometric $L_1$ adaptive control method. Section~\ref{section:sim} shows the simulation results of a quadrotor under various disturbances and uncertainties, while Section~\ref{sec:exp} presents the experimental results. Section~\ref{section: con} summarizes the work presented and provides concluding remarks.
	
	\section{Geometric Dynamics and Control of Quadrotor on SE(3) revisited}
\label{section: quad_dynamics}

\subsection{Geometric Dynamics Model}

We consider the dynamics of a quadrotor in a coordinate-free frame-work. Figure \ref{fig:quad_geo} illustrates the quadrotor with reference frames. 
The coordinate-free dynamic model is defined using a geometric representation, with the quadrotor attitude represented by a rotation matrix, $R$, in $SO(3)\coloneqq\{R\in \mathbb{R}^{3\times 3}~|~R^TR=I, det(R) = +1 \}$,  representing the rotation from body-frame to the inertial-frame. The center-of-mass position of the quadrotor, assumed to be at the geometric center of the quadrotor and denoted by $x$, is given by the vector from the inertial-frame's origin to the body-frame's origin.

The configuration space of this sytem is, $Q\coloneqq SE(3)$, with 6 degrees-of-freedom and 4 inputs corresponding to the 4 rotors. Thrust $f\in \mathbb{R}$ and moment $M\in \mathbb{R}^3$ has a one-one mapping to the 4 rotor inputs as shown in \cite{lee2011geo} and thus, can be considered as inputs to the quadrotor. 

Equations of motion for the quadrotor discussed in \cite{lee2011geo} are given below,
\begin{align}
\dot{x}&=v, \label{eq:x}\\
m\dot{v}&=mg\bm{e}_3-fR\bm{e}_3,\label{eq:dx}\\
\dot{R}& = R\hatmap{\Omega}, \label{eq:R}\\
J\dot{\Omega}& = M - (\Omega \times J\Omega), \label{eq:Omega}
\end{align}
where Table \ref{table:1} enumerates various symbols used in defining the dynamics, and the skew-symmetric  \hatmapname~$\hatmap{(\cdot)}:\mathbb{R}^3\rightarrow so(3)$ satisfies $\forall x,y\in \mathbb{R}^3, \hatmap{x}y=x\times y$. For reference, the vee-map, ${.}^{\vee}:so(3)\rightarrow \mathbb{R}^3$ is defined as the inverse of the \hatmapname~, i.e., $(\hatmap{y})^{\vee}=y$. Properties of \hatmapname~ useful in this paper are presented in Appendix A. .  

\begin{table}[t]
\centering
\begin{tabular}{l p{5cm}} 
\hline
 $m \in \mathbb{R}$ & Mass of the quadrotor \\ 
 \hline
 $J \in \mathbb{R}^{3\times 3}$ & Inertia matrix of the quadrotor with respect to the body-fixed frame \\
 \hline
$R \in SO(3) $ & Rotation matrix of the quadrotor from body-fixed frame to the inertial frame \\
\hline
$\Omega \in \mathbb{R}^3$ & Body-frame angular velocity\\
\hline
$x \in \mathbb{R}^3$ & Position vector of the quadrotor's center-of-mass in the inertial frame\\
\hline
$v \in \mathbb{R}^3$ & Velocity vector of the quadrotor's center-of-mass in the inertial frame\\
\hline
$f \in \mathbb{R}$ & Magnitude of the thrust of the quadrotor and is in the direction of $\bm{-b_3}$\\
\hline
$M\in \mathbb{R}^3$& Moment vector of the quadrotor in the body-fixed frame\\
\hline
$\bm{ e_1}, \bm{e_2}, \bm{ e_3 }\in \mathbb{R}^3$ & Unit vectors along the x,y,z directions of the inertial-frame\\
\hline
$\bm{b_1}, \bm{b_2}, \bm{b_3 }\in \mathbb{R}^3$ & Body fixed axis of the quadrotor represented in the inertial-frame; $\bm{b_3}$ is orthogonal to the plane of the quadrotor \\
\hline
$g\in \mathbb{R}$ & Acceleration due to gravity, is along the direction $\bm{e_3}$ \\
\hline
\end{tabular}
\caption{Various symbols used in representing the dynamics of the quadrotor on SE(3). }
\label{table:1}
\end{table}

\subsection{Geometric Tracking Control for Quadrotor}

The geometric tracking control presented in \cite{lee2011geo} tracks the desired quadrotor trajectory $x_d(t)$.
The position controller calculates thrust $f$ and a desired orientation $R_d$, while the attitude control calculates moment $M$ to track the desired orientation $R_d$.

The tracking errors for the attitude dynamics, with desired orientation $R_d$ and body angular velocity $\Omega_d$, are defined on the tangent bundle of $SO(3)$ as follows. The attitude tracking error is defined as,
\begin{align}
\label{eq:eRdef}
e_R = \frac{1}{2}(R_d^TR-R^TR_d)^{\vee},
\end{align}
and the angular velocity tracking error is calculated 
on $T_{R}SO(3)$ as,
\begin{align}
\label{eq:eOmdef}
e_{\Omega} = \Omega - R^TR_d\Omega_d.
\end{align}

The configuration error function between the system attitude and desired attitude is defined as,
\begin{align}
\label{eq:att_config_error}
\Psi=\frac{1}{2}Trace[ I-R_d^TR].
\end{align}
Note that, {$\Psi$ is almost globally positive definite} and upper bounded by 2.
The attitude error dynamics are then given as, see~\cite{lee2011geo} for details,
\begin{align}
\dot{e}_R &= C(R_d^TR)e_{\Omega}, \label{eq:eR1}
\end{align}
\begin{align}
J\dot{e}_{\Omega}& = J\dot{\Omega}
+ J(\hatmap{\Omega}R^TR_d\Omega_d-R^TR_d\dot{\Omega}_d), \label{eq:eOm1}
\end{align}
\begin{align}
& = M - (\Omega\times J\Omega) \nonumber\\
&\quad +J(\hatmap{\Omega}R^TR_d\Omega_d-R^TR_d\dot{\Omega}_d).\label{eq:eOm2}
\end{align}
where $C(R_d^TR)\coloneqq \frac{1}{2}(Trace[R^TR_d]I-R^TR_d)$. It can be easily noticed that choosing the control moment \highlight{$M$} as,
\begin{align}
\label{eq:Mactual}
\begin{split}
{M} = \mu + (\Omega\times J\Omega)
-J(\hatmap{\Omega}R^TR_d\Omega_d - R^TR_d\dot{\Omega}_d),
\end{split}
\end{align}
cancels the nonlinearities in \eqref{eq:Mactual}. Then a choice of $\mu$ as, 

\begin{align}
\label{eq:g_pd}
\mu = -k_Re_R - k_{\Omega}e_{\Omega},
\end{align}
for any positive constants, ${k}_R,\,{k}_{\Omega}$, would result in,
\begin{align}
J\dot{e}_{\Omega}&=-k_Re_R-k_{\Omega}e_{\Omega}.
\end{align}

For initial conditions that satisfy,
\begin{gather}
\Psi(R(0),R_d(0))<2,\\
\| e_{\Omega}(0)\|^2<\frac{2}{\lambda_{m}(J)}{k}_R(2-\Psi(R(0), R_d(0))),
\end{gather}
(where $\lambda_{m}$ is the minimum eigenvalue of the matrix and similarly, $\lambda_M$ is the maximum eigenvalue), it is shown in \cite{lee2011geo} that the zero equilibrium of the attitude tracking errors $(e_R,e_{\Omega})$ is exponentially stable. 
Moreover, a control Lyapunov candidate,
\begin{align}
V&= \frac{1}{2}e_{\Omega}\cdot Je_{\Omega} + {k}_R\Psi(R,R_d) + c_2e_R\cdot e_{\Omega},
\end{align}
is defined to show that the above error dynamics is exponentially stable for
\begin{equation*}
c_2<\bigg\{k_{\Om},\frac{4k_{\Omega}k_R\lambda_m(J^2)}{k^2_{\Om}\lambda_M(J)+4k_R\lambda_m(J^2)},\sqrt{k_R\lambda_m(J)} \bigg\},
\end{equation*}
by proving that $\dot{V}\leq \eta^TW\eta$, where $\eta=[e_R \quad e_{\Omega}]^T$ and $W$ is a positive definite matrix. A detailed proof can be found in \cite{lee2011geo}. 

Having discussed the geometric attitude control, the following section describes the attitude dynamics of quadrotor in the presence of disturbances and model uncertainties.

	\section{Effects of Model Uncertainties and Disturbances on the Attitude Dynamics}
\label{sec:disturb}
As discussed in the previous section, the control moment $M$ in \eqref{eq:Mactual} ensures that the zero equilibrium of the error dynamics in \eqref{eq:eR1}-\eqref{eq:eOm2} is exponentially stable. This controller assumes an accurate dynamical model of the quadrotor. However, presence of any uncertainties in the model properties, like in mass, $m$, and in inertia, $J$, can result in large tracking errors and potential instability. External disturbances on the system can also result in similar adverse effects. 

The unknown external disturbances can be captured in the attitude dynamics of the quadrotor \eqref{eq:R}-\eqref{eq:Omega}. In particular, the attitude dynamics along with the external disturbances are defined as,
\begin{align}
\dot{R} & = R\hatmap{\Omega}, \label{eq:R_real}\\
\Jtrue\dot{\Omega} &= M - \Omega\times \Jtrue\Omega + \theta_e, \label{eq:Omega_real}
\end{align}
where $\theta_e$ represents the unknown external disturbance
and $\Jtrue$ is the true (unknown) inertia of the quadrotor. The corresponding attitude error dynamics \eqref{eq:eR1}, \eqref{eq:eOm2}, whose errors are defined as \eqref{eq:eRdef} and \eqref{eq:eOmdef}, can be modified and represented as below,
\begin{align}
\dot{e}_R &= C(R_d^TR)e_{\Omega}, \label{eq:eR_real}\\
\Jtrue\dot{e}_{\Omega}& = M - (\Omega\times \Jtrue\Omega) \nonumber\\
&\quad +\Jtrue(\hatmap{\Omega}R^TR_d\Omega_d-R^TR_d\dot{\Omega}_d) + \theta_e. \label{eq:eOm_real_2}
\end{align}

\begin{remark}
Due to the above structure, the external disturbance $\theta_e$ becomes a matched uncertainty~\cite[Chapter~2]{hovakimyan2010l1}.
\end{remark}

Additionally, the control moment \highlight{$M$} in \eqref{eq:Mactual} assumes accurate knowledge of inertia of the quadrotor. However, if the true inertia, $\Jtrue$, of the quadrotor is not same as the nominal inertia of the quadrotor, $\Jnom$, it results in further uncertainties in the closed-loop system.

In particular, substituting for \highlight{$M$} from \eqref{eq:Mactual} in \eqref{eq:eOm_real_2}, we obtain,
\begin{align}
\Jtrue\dot{e}_{\Omega}&=\mu + \theta_e \nonumber \\
&+ {\scriptsize\underbrace{[-\Omega\times(\Jtrue-\Jnom)\Omega + (\Jtrue-\Jnom)(\hatmap{\Omega}R^TR_d\Omega_d-R^TR_d\dot{\Omega}_d)]}_{\triangleq \theta_m}}, \label{eq:Om_real_3}
\end{align}
where, $\theta_m$ is zero when $\Jnom = {\Jtrue}$.
\highlight{We define $\delta J$ relating the true inertia ${\Jtrue}$ and the nominal inertia $\Jnom$ as,
$\delta J := \Jnom\Jtrue^{-1} - I$.}

Thus, \eqref{eq:Om_real_3} can be written in the following manner,
\begin{align}
\Jnom\dot{e}_{\Omega}&= \mu + \underbrace{(\delta J)\mu + \Jnom\Jtrue^{-1}(\theta_m+\theta_e)}_{\triangleq \theta} \label{eq:theta_def}
\end{align}
where $\theta$ is the combined expression for the uncertainties and disturbances. Therefore, the closed-loop attitude error dynamics of the quadrotor with control moment $M$ defined in \eqref{eq:Mactual}, along with the model uncertainties and disturbances can be given as,
\begin{subequations}
\label{eq:e_real_final}
\begin{empheq}[box=\widefbox]{align}
  \dot{e}_R&=C(R_d^TR)e_{\Omega}, \label{eq:eR_real_final}\\
  \Jnom\dot{e}_{\Omega}&= \mu + \theta. \label{eq:eOm_real_final}
\end{empheq}
\end{subequations}

For $\theta=0$ in \eqref{eq:eOm_real_final}, representing no model uncertainity, a choice of $\mu$ as a PD control similar to \eqref{eq:g_pd}, would result in a stabilizing control. However, if $\theta\neq 0$, the choice of PD control alone will not be sufficient to guarantee stability. A choice of $\mu$ that can cancel the uncertainty $\theta$ would be helpful in achieving stability. However, {$\theta$ is unknown and this cannot be done}. A similar approach for biped robots is developed in \cite{quan2015bi}, however, they do not address the case of dynamics evolving on manifolds. Also note that the uncertainty is a nonlinear function of states $R,\Omega$ and control input $\mu$. In the next section, we propose a Geometric $L_1$ adaptation to predict the uncertainty $\theta$.

	\section{Geometric $L_1$ Adaption for Attitude Tracking Control}
\label{section: control}


In section \ref{sec:disturb}, the attitude dynamics were described in the presence of model uncertainties. In particular, \eqref{eq:eR_real_final} \& \eqref{eq:eOm_real_final} present the attitude error dynamics along with the uncertainty $\theta$ and input $\mu$. (Note that the actual control moment $M$ is calculated using \eqref{eq:Mactual}.) In this section, we proceed to present a geometric $L_1$ adaptation law to estimate the uncertainty $\theta$, and compute the input $\mu$, to track a desired time-varying trajectory $(R_d,\Omega_d)$. 

\begin{remark}
The dynamical model given in \eqref{eq:eOm_real_2} with the controller in \eqref{eq:Mactual} resulting in the closed-loop system in \eqref{eq:e_real_final} is referred to as \qmodel.
States and inputs corresponding to the \textit{\qmodel} are given below,
\begin{equation*}
\boxed{R,\Omega, e_R, e_{\Omega}, \theta, M.}
\end{equation*}
\end{remark}

\subsection{\Refsys}

$L_1$ control architecture employs a \refsys~(also referred as reference/nominal system or state predictor) to predict the uncertainty $\theta$ in the system. In this paper, we consider a \refsys\, with a nominal inertia matrix, $\Jnom$, and without any disturbances $\theta_e$. Let $\hat{R}$ be the attitude of the \refsys~and $\hat{\Omega}$ be the body-angular velocity of the \refsys. Dynamics of the \refsys\, is written as,
\begin{align}
\dot{\hat{R}}&=\hat{R}\hatmap{\hat{\Omega}}, \label{eq:R_ref}\\
\Jnom\dot{\hat{\Omega}}&=\hat{M} - \hat{\Omega}\times \Jnom\hat{\Omega}, \label{eq:Om_ref}
\end{align}
with control moment $\hat{M}$ defined similar to \eqref{eq:Mactual}, i.e.,
\begin{align}
\hat{M}&=\hat{\mu} + (\hat{\Omega}\times \Jnom\hat{\Omega})\nonumber\\
&\quad - \Jnom(\hatmap{\hat{\Omega}}\hat{R}^TR_d\Omega_d - \hat{R}^TR_d\dot{\Omega}_d). \label{eq:Mhat}
\end{align}

Here, $R_d(t)$, $\Omega_d(t)$ are the same desired trajectory (attitude and body-angular velocity) considered in section \ref{sec:disturb}. 
Similar to the attitude error vectors defined for the \qmodel, we define configuration errors for the \refsys. The attitude tracking error is defined as,
\begin{align}
\hat{e}_R&=\frac{1}{2}(R_{d}^{T}\hat{R}-\hat{R}^{T}R_{d})^{\vee}, \label{eq:eRhat}
\end{align}
and the angular velocity tracking error on $T_{\hat{R}}SO(3)$ is,
\begin{align}
\label{eq:eOmhat}
\hat{e}_{\Omega}&= \hat{\Omega}-\hat{R}^TR_d\Omega_d.
\end{align}

The error dynamics for the  \refsys\, is similar to \eqref{eq:eR1}, \eqref{eq:eOm2}. We present the attitude error dynamics with control moment defined in \eqref{eq:Mhat} below,

\begin{subequations}
\begin{empheq}[box=\widefbox]{align}
  \dot{\hat{e}}_R&=C(R_d^T\hat{R})\hat{e}_{\Omega}, \label{eq:deRhat}\\
  J\dot{\hat{e}}_{\Omega}&=\hat{\mu}   \label{eq:deOmhat}
\end{empheq} \label{eq:tildeerrordyn}
\end{subequations}

Comparing \eqref{eq:deOmhat} to \eqref{eq:eOm_real_final}, we notice the presence of additional term $\theta$ (uncertainty) in the \qmodel. $L_1$ adaptation is used to estimate this uncertainty, with $\hat{\theta}$ denoting the uncertainty. Discussion regarding the uncertainty estimation and its relevance to the control is presented below,
{ 
\begin{enumerate}[{\it (i)}]
\item $L_1$ adaptation makes use of both the \qmodel~ and the \refsys~ to estimate the uncertainty $\hat{\theta}$. 
\item In-order for the \refsys~ to account for the uncertainty in the \qmodel, estimated uncertainty $\hat{\theta}$ (or a transported version of $\hat{\theta}$ in case of manifolds) is included in the control input $\hat{\mu}$.    
\item Control inputs $\mu$ (\qmodel) and $\hat{\mu}$ (\refsys) are used to track a desired trajectory while canceling the uncertainty. 
\item Uncertainty is countered by including $(-\hat{\theta})$ in the control inputs $\mu$ and $\hat{\mu}$. However, $\hat{\theta}$ typically contains high frequency components due to fast estimation. $L_1$ adaptive control architecture is used to decouple estimation and adaption \cite{Cao2006adapt}.  A low-pass filter is used to exclude the high frequency content in the input. Thus, $(-C(s)\hat{\theta})$ is included in the control inputs, where $C(s)$ is a low pass filter and with \highlight{$\| C(0)\| = 1$}. The low-pass filter is key to the trade-off between the performance and robustness. 

\end{enumerate}}

\begin{remark}
The \refsys~is distinguished using the superscript $\hat{\cdot}$. Therefore, states and inputs corresponding to the \refsys\, are as given below,
\begin{equation*}
\boxed{\hat{R},\hat{\Omega}, \hat{e}_R, \hat{e}_{\Omega}, \hat{\theta}, \hat{M}.}
\end{equation*}
Table \ref{table:2} presents the different notations used in this paper.
\end{remark}

\begin{table}[t]
\centering
\begin{tabular}{m{2cm}  c  m{4cm}} 
\hline 
Symbol & Example	  & Model/Errors \\
\hline \hline 
no sub/super-script & $R$  & \Qmodel \\ [1ex]
\hline
${.}_d$ & ${R}_d$  & Desired Trajectory \\ [1ex]
\hline
$\widehat{.}$ & $\widehat{R}$  & \Refsys \\ [1ex]
\hline \hline
no sub/super-script & ${e}_R$  & Error between \Qmodel~\& Desired Trajectory  \\ [1ex] 
\hline 
$\widehat{.}$ & $\widehat{e}_R$  & Error between \Refsys ~\& Desired Trajectory  \\ [1ex] 
\hline 
$\widetilde{.}$ & $\widetilde{e}_R$  & Error between \Qmodel~\& \Refsys \\ [1ex] 
\hline 
\hline
& $\Jtrue$ & True inertia \\
\hline
& $\Jnom$ & Nominal inertia used by the control law\\
\hline 
\end{tabular}
\caption{List of notations used to represent various models and errors in this paper.}
\label{table:2}
\end{table}

\subsection{Errors between the \Qmodel\, and the \Refsys\,}

In the previous sections, we presented the 
\qmodel~and the corresponding \refsys~to mimic the actual system and its uncertainties. The goal of our work is to present an adaptive controller to reduce the differences between the \qmodel~and the \refsys. We, also present a tracking controller for the \refsys, that enables the \qmodel~ to track the desired trajectory.

Attitude tracking errors between the \qmodel\, and the desired time-varying trajectory are $e_R$ and $e_{\Omega}$. Errors between the \refsys\, and the desired trajectory are $\hat{e}_R$ and $\hat{e}_{\Omega}$. 
Similarly, we define a new set of attitude errors to capture the difference between the \qmodel\, and the \refsys~as below, 
\begin{align}
\tilde{e}_R &=\frac{1}{2}(R^T\hat{R}-\hat{R}^TR)^{\vee}, \label{eq:eRtilde}\\
\tilde{e}_{\Omega}&=\hat{\Omega}-\hat{R}^TR\Omega, \label{eq:eOmtilde}
\end{align}
and the corresponding configuration error function is given as,
\begin{align}
\tilde{\Psi}&=\frac{1}{2}Trace[I-R^T\hat{R}].\label{eq:Psitilde}
\end{align}

{From \eqref{eq:eOm_real_final} and \eqref{eq:Mactual}, the uncertainty $\theta$ can be canceled out (if it was known) by the control moment.} This shows that $\theta$ lies in the same dual space as the moment \highlight{$M$}, i.e., $\theta \in T^*_RSO(3)$, and similary, the predicted uncertainty is in the dual space $T^*_{\hat{R}}SO(3)$. The difference between the actual uncertainty and the predicted uncertainty is thus calculated by transporting ${\theta}$ to the space of $\hat{\theta}$,
\begin{align}
\tilde{\theta}&=\hat{\theta} - \hat{R}^TR\theta. \label{eq:thetatilde}
\end{align}

Having presented the error definitions for orientation $\tilde{e}_R$, angular velocity $\tilde{e}_{\Omega}$ and uncertainty $\tilde{\theta}$ between the \qmodel~and the \refsys, we now present the control design for $\mu$ and $\hat{\mu}_1$ for the attitude dynamics of \qmodel, and \refsys~respectively, such that these errors, $(\tilde{e}_R,\tilde{e}_{\Omega})$, exponentially reach {an arbitrarily  small neighborhood of the origin (0,0)}.

\subsection{Geometric Attitude Tracking Control for quadrotor with uncertainty} 

{In this section, we define the control inputs $\mu$ and $\hat{\mu}$. Intuition behind the definition of the control inputs is presented below,
\begin{itemize}
\item \Refsys:

\begin{enumerate}[{\it (i).}]
\item Control moment $(\hat{M})$ and the control input $(\hat{\mu})$ for the \refsys~is given in \eqref{eq:controlMomentMhat}
\item  $\hat{\mu}$ is designed to track the desired trajectory (through $\hat{\mu}_1$), i.e.,  $(\hat{R},\hat{\Omega})\rightarrow (R_d,\Omega_d)$, alternately $(\hat{e}_{R},\hat{e}_{\Omega})\rightarrow (0,0)$
\item  $\hat{\mu}$ also includes the terms to account for the true uncertainty and low-pass filtered estimate of the uncertainty to counter the uncertainty (through $\hat\mu_2$)
\end{enumerate}
\end{itemize}}

{
\begin{itemize}
\item \Qmodel:
\begin{enumerate}[{\it (i).}]
\item Control moment $({M})$ and the control input $({\mu})$ for the \qmodel~is given in \eqref{eq:controlMomentM_first}
\item $\mu$ is defined such that the \qmodel~tracks the \refsys, i.e., $(R,\Omega)\rightarrow (\hat{R}, \hat{\Omega})$, alternately $(\tilde{e}_R, \tilde{e}_{\Omega})\rightarrow(0,0)$
\item Low-pass filtered uncertainty estimate, $-C(s)\hat{\theta}$,  is included to cancel the uncertainty in the \qmodel. However, $-C(s)\hat{\theta}$ has to be transported to the ${\mu}$ space. 
\item $\mu$ is chosen such that the error dynamics for the errors in \eqref{eq:eRtilde}, \eqref{eq:eOmtilde} are feedback linearized and results in $J\dot{\tilde{e}}_{\Omega}=-\tilde{k}_R\tilde{e}_R-\tilde{k}_{\Omega}\tilde{e}_{\Omega}+\Rho\tilde{\theta}$. 
\end{enumerate}
\end{itemize}
}

\begin{tcolorbox}[width=\columnwidth,colback={LightGray},title={Control moment, $(\hat{M})$, for the \refsys},colbacktitle=white,coltitle=black]
\begin{subequations}
\begin{align}
\hat{M}&=\hat{\mu} + (\hat{\Omega}\times J\hat{\Omega})\nonumber\\
&\quad - J(\hatmap{\hat{\Omega}}\hat{R}^TR_d\Omega_d - \hat{R}^TR_d\dot{\Omega}_d),\\
 \hat{\mu}&= \hat{\mu}_1+\hat{\mu}_2,\\
 \hat{\mu}_1 &= -\hat{k}_R\hat{e}_R - \tilde{k}_{\Om}\hat{e}_{\Om},\label{hatmu1} \\
 \hat{\mu}_2 &=  J\hat{R}^TRJ^{-1}R^T\hat{R}\hat{\theta}-C(s)\hat{\theta}.
 \end{align}  \label{eq:controlMomentMhat}
\end{subequations}
\end{tcolorbox}

\begin{tcolorbox}[width=\columnwidth,colback={LightGray},title={Control moment, $(M)$, for the \qmodel},colbacktitle=white,coltitle=black]
\begin{subequations}
\begin{align}
  M& = \mu +(\Omega\times J\Omega) \nonumber \\
	&\quad \quad -J(\hatmap{\Omega}R^TR_d\Omega_d - R^TR_d\dot{\Omega}_d),\\
 \mu&= \mu_1+\mu_2, \\
 \mu_1 & = JR^T\hat{R}J^{-1}(\hat{\mu}_1{-}C(s)\hat{\theta}+\tilde{k}\tilde{e}_R{ + }\tilde{k}_{\Omega}\tilde{e}_{\Omega}), \label{mu1}\\
\mu_2 & =  JR^T\hat{R}\hatmap{\tilde{e}}_{\Om}\hat{R}^TRe_{\Om}.
 \end{align}  \label{eq:controlMomentM_first}
\end{subequations}
\end{tcolorbox}

In \eqref{eq:controlMomentMhat}, \eqref{eq:controlMomentM_first} $\hat{k}_R$, $\tilde{k}_R$ and $\tilde{k}_{\Omega}$ are positive constants. {The resulting attitude error dynamics  are presented in Appendix B.} 
 
\begin{remark}
The relation between $\mu$ and $\hat{\mu}$ can be shown as follows,
\begin{align}
\mu & = JR^T\hat{R}J^{-1}(\hat{\mu}+\tilde{k}_R\tilde{e}_R+\tilde{k}_{\Omega}\tilde{e}_{\Omega})\nonumber\\
&\quad\quad -R^T\hat{R}\hat{\theta} + JR^T\hat{R}\hatmap{\tilde{e}}_{\Omega}\hat{R}^TRe_{\Omega}. \label{eq:controlMomentM}
\end{align}
This relation is used to calculate $\dot{\tilde{e}}_{\Omega}$ in Appendix B. \highlight{Note that the presence of $\tilde{k}_{\Omega}$ in both $\hat{\mu}_1$ (34c) and $\mu_1$ (35c), ensures that the relation in \eqref{eq:controlMomentM} is obtained after simplification.}
\end{remark} 

\begin{figure}
\centering
\includegraphics[width=\columnwidth]{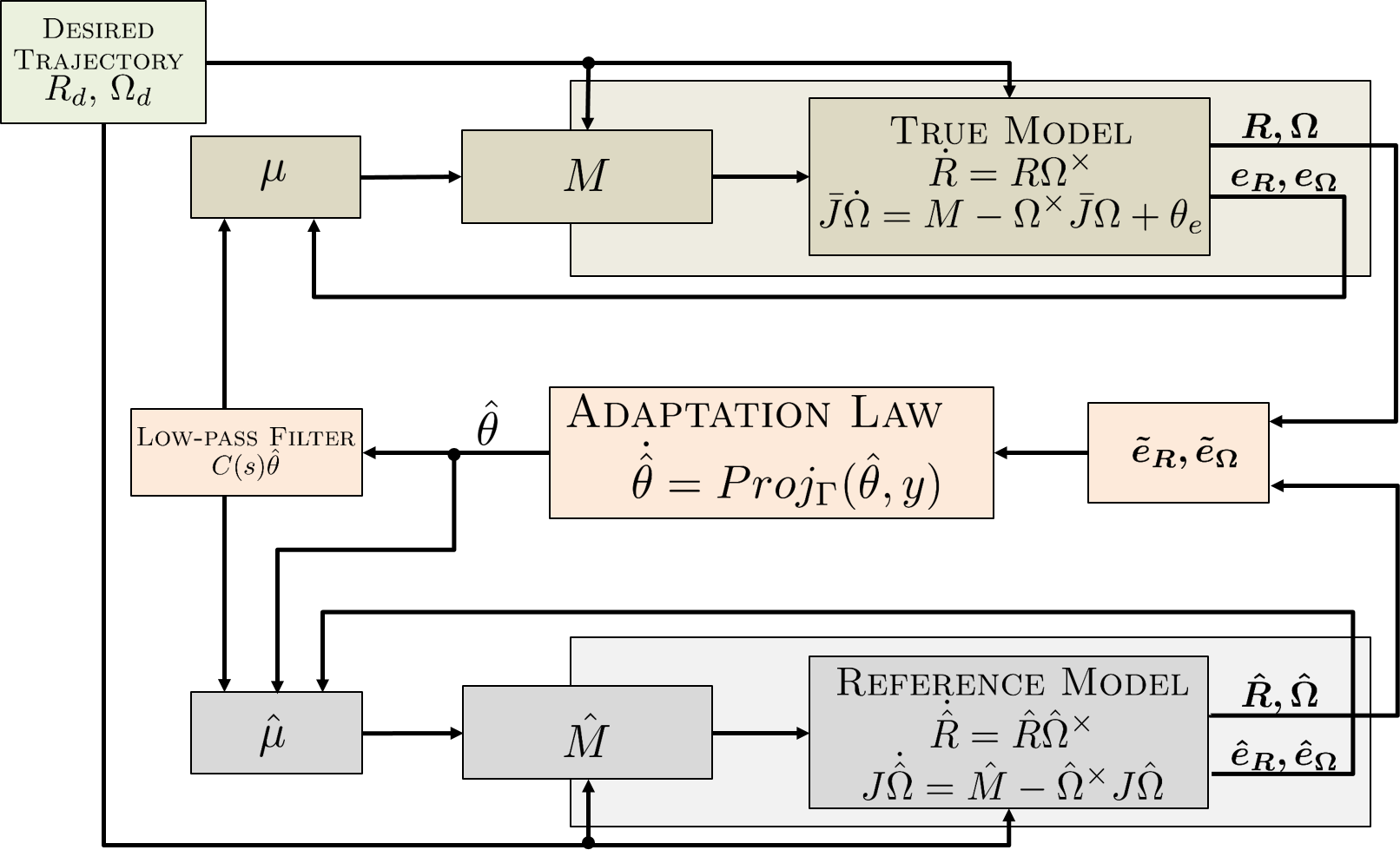}
\caption{$L_1$ Adaptive scheme on the Geometric Attitude Control. The \qmodel~captures the dynamics of the quadrotor and represents the physical plant. A $4^{th}$ order Range-Kutta method is used to simulate the dynamics of the \Qmodel~ and \Refsys {\it s} in the Numerical Simulations and Euler integration is used to simulate the Reference model in the experiments.}
\label{fig:l1scheme}
\end{figure}

\begin{thm}
\label{prop:1}
Consider the control moments $M$ and $\hat{M}$ defined in \eqref{eq:controlMomentM_first} and \eqref{eq:controlMomentMhat} for the \qmodel\,and the \refsys\, respectively and 
\begin{enumerate}[{\it (i)}]
\item the adaptation law using $\Gamma-$Projection as,
\begin{equation}
\label{eq:adapt}
\dot{\hat{\theta}} = Proj_{\Gamma}(\hat{\theta},y),
\end{equation}
with the projection operator as defined in Remark~\ref{remark:gamma},
\item the definition of $y$ given as,
\begin{align}
\label{eq:y}
y = -(\Rho^T\tilde{e}_{\Omega}+c\Rho^TJ^{-T}\tilde{e}_R),
\end{align} where 
\begin{align}
\label{eq:Rho_def}
\Rho = J\hat{R}^TRJ^{-1}R^T\hat{R},
\end{align} and
\item with the initial condition that satisfies,
\begin{gather}
\T{\Psi}(\hat{R}(0),R(0))<2,\\
\| \T{e}_{\Omega}(0)\|^2<\frac{2}{\lambda_{min}(J)}\T{k}_R(2-\T{\Psi}(\hat{R}(0)R(0))), \label{eq:eOmtilde_init}
\end{gather}
\end{enumerate}
then the attitude tracking error $(\tilde{e}_R, \tilde{e}_{\Omega})$, defined in \eqref{eq:eRtilde}, \eqref{eq:eOmtilde}, is \textit{exponential input-to-state stable (e-ISS)} \cite{kolathaya2018input} in the sense of Lyapunov. 
\begin{proof}
Proof is given in Appendix C. 
\end{proof}
\end{thm}
\begin{remark}
\label{remark:gamma}
Definition of the $\Gamma-$Projection operation is given in \cite{lavretsky2011projection} as,
  \begin{equation}
    Proj_{\Gamma}(\hat{\theta},y)= 
    \begin{cases*}
       \Gamma y - \Gamma\frac{\bigtriangledown f(\hat{\theta})(\bigtriangledown f(\hat{\theta}))^T}{(\bigtriangledown f(\hat{\theta}))^T(\bigtriangledown f(\hat{\theta}))}\Gamma y f(\hat{\theta}), & \\
      \quad\quad \text{if } f(\hat{\theta}) > 0 \wedge y^T\Gamma \bigtriangledown f(\hat{\theta})> 0, & \\
      \Gamma y, \quad\quad \text{otherwise},
    \end{cases*}
  \end{equation}
where, $f(\hat{\theta}):\mathbb{R}^3\rightarrow \mathbb{R}$ is any convex function.
\end{remark}
Figure \ref{fig:l1scheme} presents the control architecture with $L_1$ adaptation. Also, note the presence of low-pass filter after the adaptation law. This is integral to the $L_1$ adaptive controller, where the high frequency noise is filtered from the system input~\cite{cao2006design}.  The $L_1$  adaptive controller consists of the reference model \eqref{eq:tildeerrordyn}, adaptation law \eqref{eq:adapt}, and the control law \eqref{eq:controlMomentM_first}. 

\subsection{Position Control}
Now that we have an attitude controller, we can use the position controller from \cite{lee2011geo}.
For some smooth position tracking command $x_{d}(t) \in \mathbb{R}^{3}$, we can define the position and velocity tracking errors as,
\begin{equation}
\label{eq:ex}
e_{x} = x - x_{d},
\end{equation}
\begin{equation}
\label{eq:ev}
e_{v} = v - v_{d}.
\end{equation}

The desired thrust vector for the quadrotor is computed as
\begin{equation}
\label{eq:thrust_vector}
\vec{F} = -k_{x}e_{x}-k_{v}e_{v}-mg\bm{e}_{3}+m\ddot{x}_{d}
\end{equation}

The quadrotor thrust magnitude is then given by,
\begin{equation}
\label{eq:f}
f =- \vec{F} \cdot R\bm{e}_{3},
\end{equation}
where $k_{x}$ and $k_{v}$ are positive constants and the thrust is in direction $-\bm{b_{3}}$.
The desired orientation and angular velocity are given by,
\begin{equation}
R_{d}=\begin{bmatrix}
\bm{b_{1d}} &~~~\bm{b_{3d}} \times \bm{b_{1d}}&~~~\bm{b_{3d}}
\end{bmatrix},
\quad
\Omega_{d}^{\times} = R_{d}^{T}\dot{R}_{d},
\end{equation}
where $\bm{b_{3d}}$ is chosen as,
\begin{equation}
\label{eq:b3d}
\bm{b_{3d}}=- \frac{-k_{x}e_{x}-k_{v}e_{v}-mg\bm{e}_{3}+m\ddot{x}_{d}}{\|-k_{x}e_{x}-k_{v}e_{v}-mg\bm{e}_{3}+m\ddot{x}_{d}\|},
\end{equation}
and is selected to be orthogonal to $\bm{b_{1d}}$, such that $R_{d} \in SO(3)$ - see \cite{lee2011geo}, \highlight{\cite[chapter 11]{quan2017introduction}} for more details.

\begin{remark}
A similar $L_1$ approach, used in the attitude control, can also be used to deal with disturbances and uncertainties in position control as studied in \cite{zuo2014adapt}. Since position and velocity evolve in the Cartesian space, we use a traditional (non-geometric) position $L_1$ adaptive control in the simulations and experiments discussed in the later sections.  
\end{remark}

%
\section{Numerical Validation}
\label{section:sim}

In this section, we present numerical examples to validate the controller presented in the previous section. We discuss the performance of geometric $L_1$ control compared to geometric control \textit{without} $L_1$ \highlight{(geometric PD control) \cite{lee2011geo}}. Two different scenarios are considered to evaluate the performance of geometric $L_1$ controller. {We also compare geometric $L_1$ control performance to Euler angle $L_1$ control to present a better insight into the performance of geometric control.}


The following system properties are considered, $m = 1.129kg$, $J =diag([6.968,\,6.211,\,10.34])\times10^{-3} kg \cdot m^{2}$,
and the following control parameters are used,
$ k_{x} = 4, k_{v} = 3.2, k_{R}=\tilde{k}_R = 2, k_{\Omega}=\hat{k}_{\Omega}=\tilde{k}_{\Omega} = 0.25.
$
A first-order low-pass filter of form, $
C(s) = \frac{a}{s + a}$, is used with $a=2$ and an adaptation gain of $\Gamma = 10^{6}I_{3\times3}$.

We present simulations for a quadrotor tracking a desired circular trajectory given below,
\begin{equation}
\label{eq:circ_traj}
\highlight{x(t) = [\rho\cos(\omega t),\, \rho\sin(\omega t),\, 0]^T,\, \psi(t)=0,} 
\end{equation}

where $\highlight{\rho=1}$, $\omega = 2$. The two different scenarios considered are, (i) a constant external disturbance and (ii) a time varying disturbance. In both scenarios, we consider a mass $m_a=0.5kg$ attached to the quadrotor at $r=[0.2, 0.2, 0.2]^Tm$ in the body frame, this added mass and its location is unknown to the controller. Presence of this added mass will also result in a moment about the center-of-mass due to gravity. Thus the external disturbance due to the mass is given as,
\begin{align}
\label{eq:theta_e_1}
\theta_e = m_ag(r\times R^T\bm{e_3})
\end{align}
and the inertia is
\rebuttal{
\begin{align}
{\Jtrue}=J+J_{m_a}, \label{eq:true_inertial_calc}
\end{align}
where, $\Jtrue$ is the true inertia used in simulating the dynamics of the \qmodel\, and $\Jnom$ is the nominal inertia used to calculate the control moment and $J_{m_a}= -m_a(\hatmap{r})^2$ is inertia tensor due to $m_a$ in the body-frame.}

\begin{figure}
    \centering
    \begin{subfigure}[t]{0.5\columnwidth}
        \centering
        \includegraphics[width=\columnwidth]{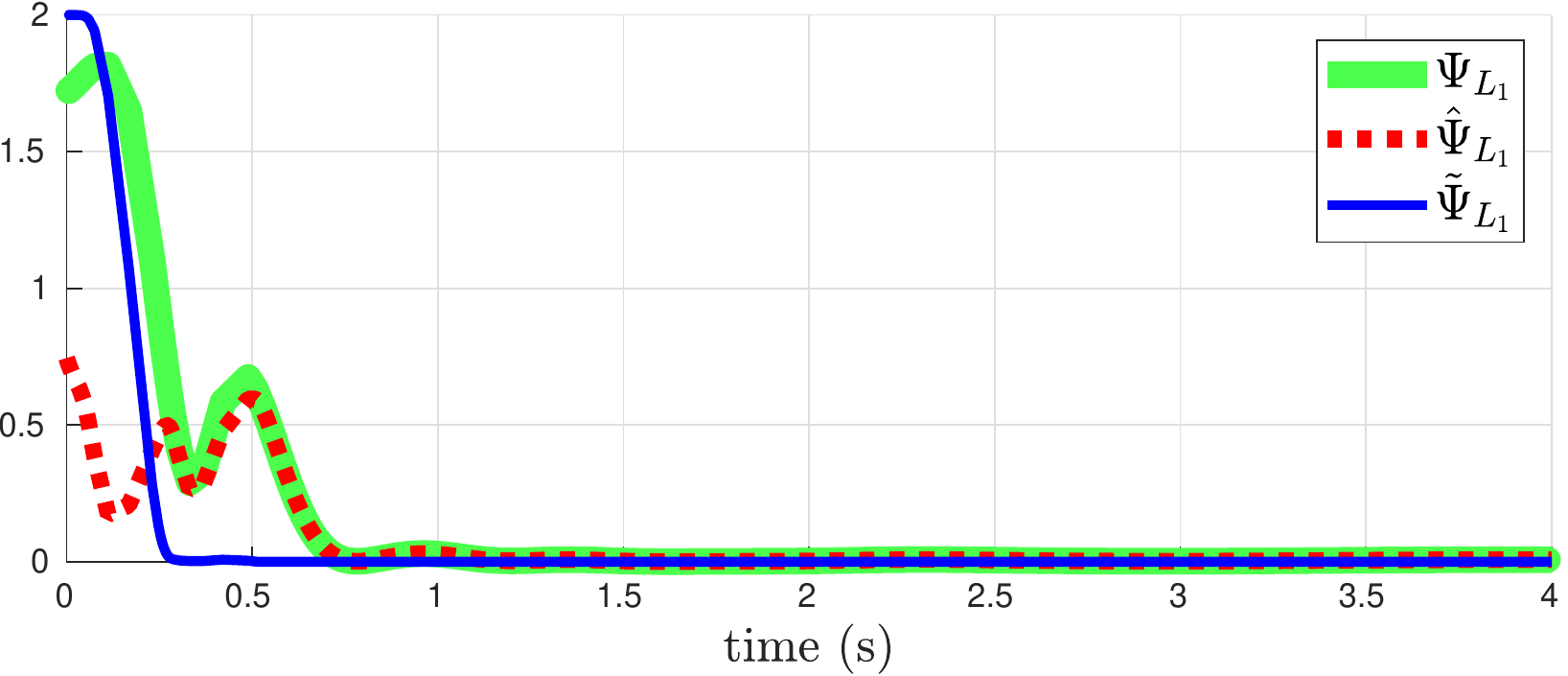}
        \caption{Configuration Errors for $L_1$: $\Psi\text{ vs }\hat{\Psi}\text{ vs }\tilde{\Psi}$}
        \label{fig_sim_1a}
    \end{subfigure}%
    ~ 
    \begin{subfigure}[t]{0.5\columnwidth}
        \centering
        \includegraphics[width=\columnwidth]{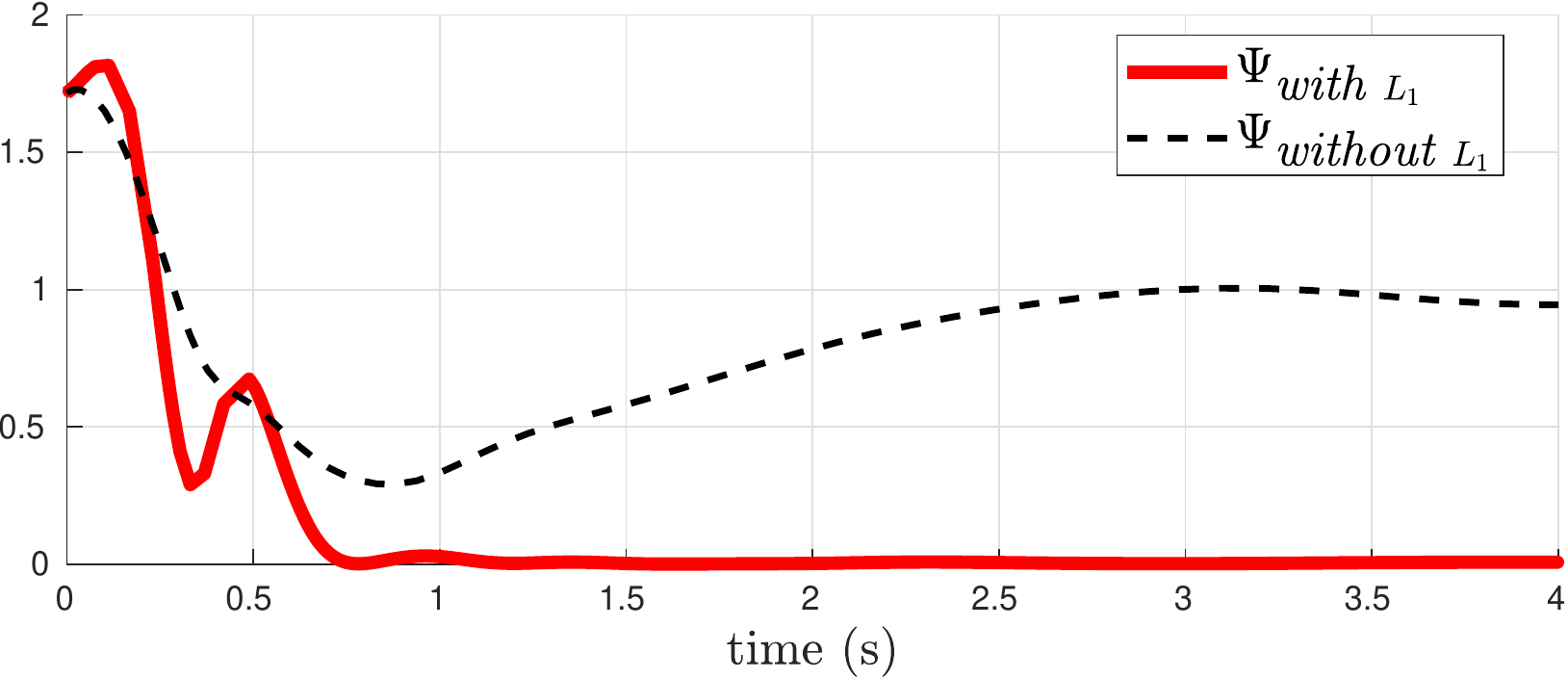}
        \caption{Configuration Errors: $L_1$ vs without $L_1$}
        \label{fig_sim_1b}
    \end{subfigure}
    
    \begin{subfigure}[t]{0.5\columnwidth}
        \centering
        \includegraphics[width=\columnwidth]{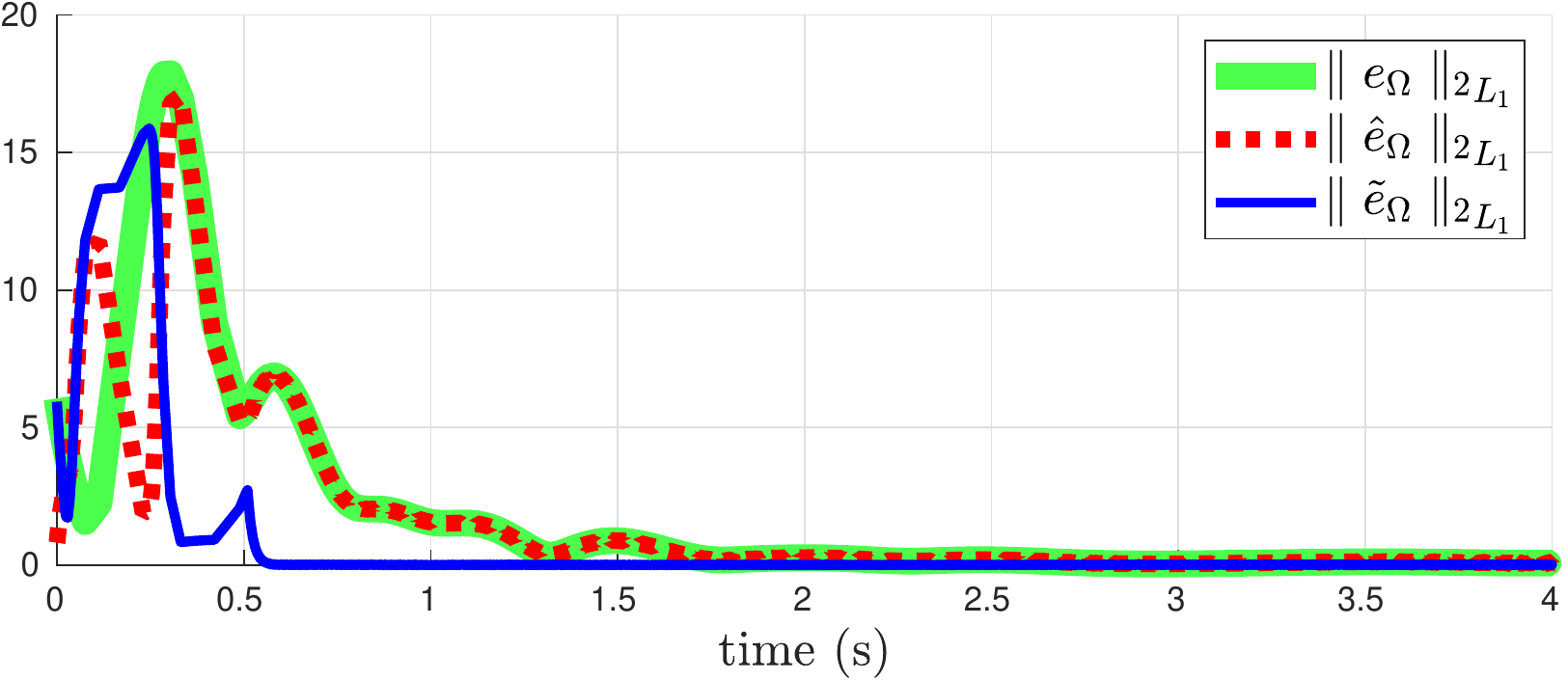}
        \caption{Angular Velocity Errors for $L_1$: $e_{\Omega}\text{ vs }\hat{e}_{\Omega}\text{ vs }\tilde{e}_{\Omega}$}
        \label{fig_sim_1c}
    \end{subfigure}%
    ~ 
    \begin{subfigure}[t]{0.5\columnwidth}
        \centering
        \includegraphics[width=\columnwidth]{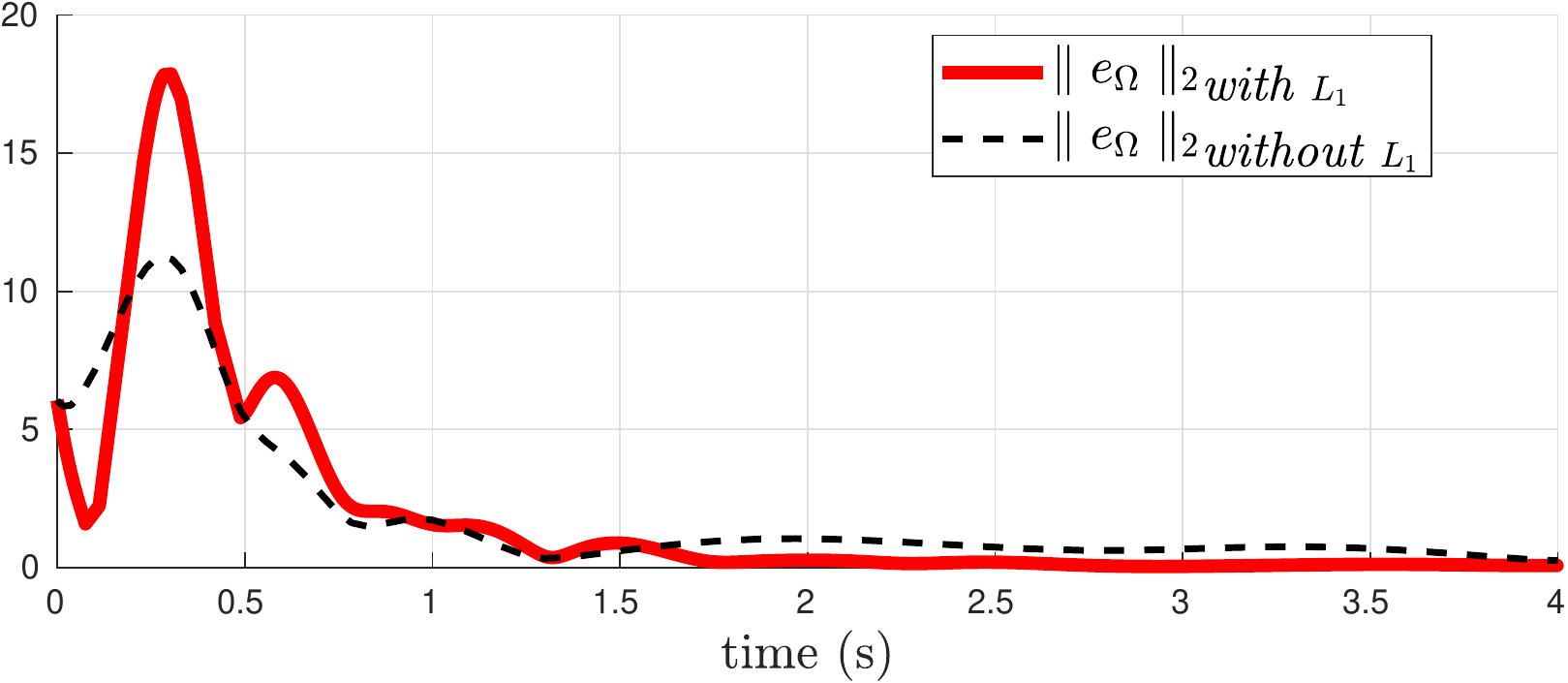}
        \caption{Angular Velocity Errors: $L_1$ vs without $L_1$}        
        \label{fig_sim_1d}
    \end{subfigure}
    \caption{\textit{Circular Trajectory - Case I}: Comparison of errors between \qmodel\, and \refsys\, in the presence of constant external disturbance $\theta_t = [.95,.25,-.5]^T$ and model uncertainty ($m_a=0.5kg$ at $r = [0.2,0.2,0.2]\highlight{^T}m$). Subfigures \ref{fig_sim_1a} \& \ref{fig_sim_1c} show that $\tilde{\Psi}$ and $\|\tilde{e}_{\Omega}\|$ decrease to zero even though the \qmodel~and \refsys~are initialized to different values (see \ref{fig_sim_1a}). Subfigures \ref{fig_sim_1b} and \ref{fig_sim_1d} show that the errors do not converge to zero in case of geometric  control \textit{without} $L_1$ (i.e, Geometric control without the $L_1$ adaptation law).}
    \label{fig:circ_case_b_errors}
\end{figure}

\begin{figure}
\centering
\includegraphics[width=1\columnwidth]{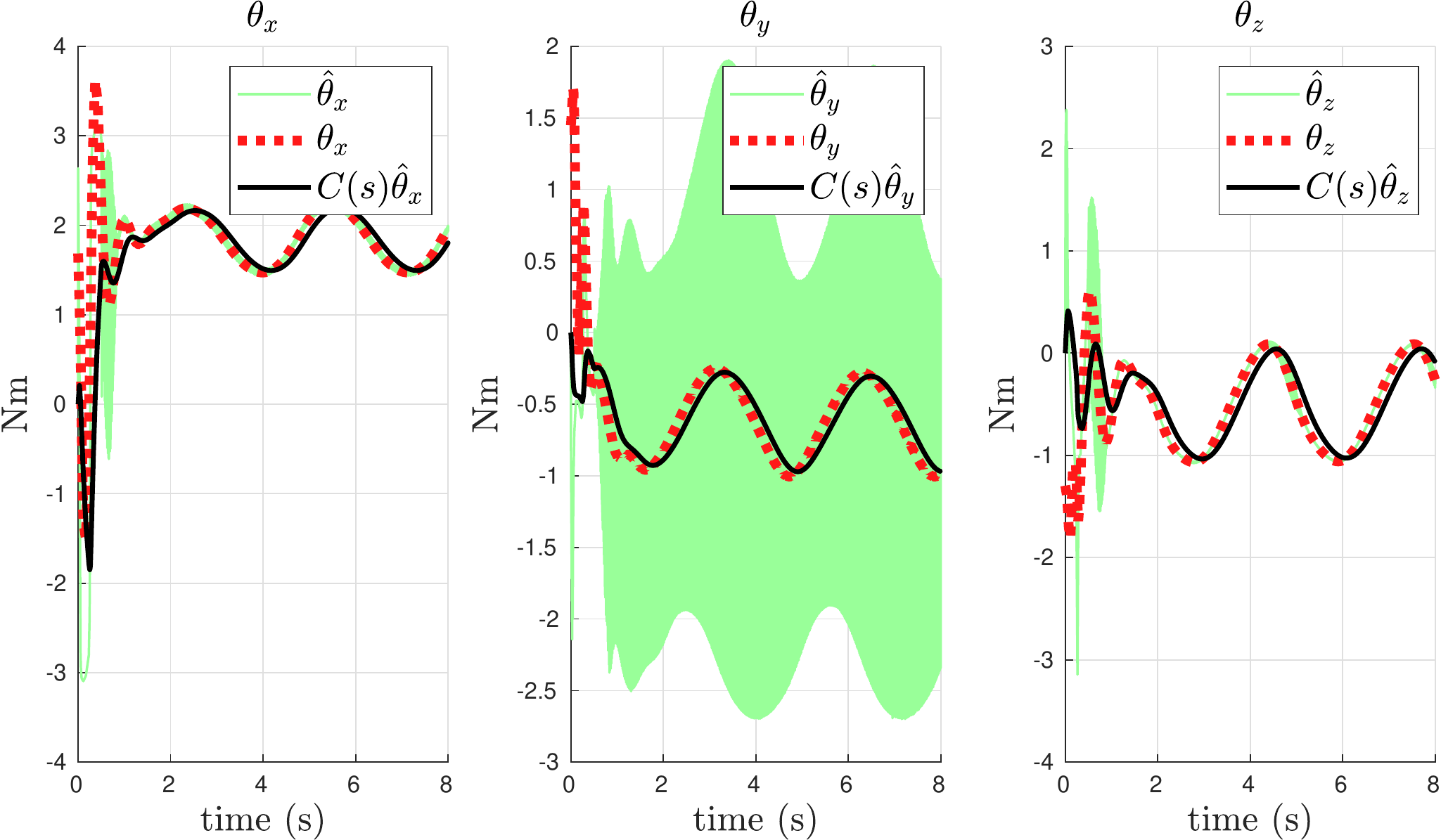}
\caption{\textit{Circular Trajectory - Case I}: Uncertainty in the system calculated using \eqref{eq:theta_def}. }
\label{fig:circ_case_b_theta}
\end{figure} 

\begin{figure}
\centering
\includegraphics[width=1\columnwidth]{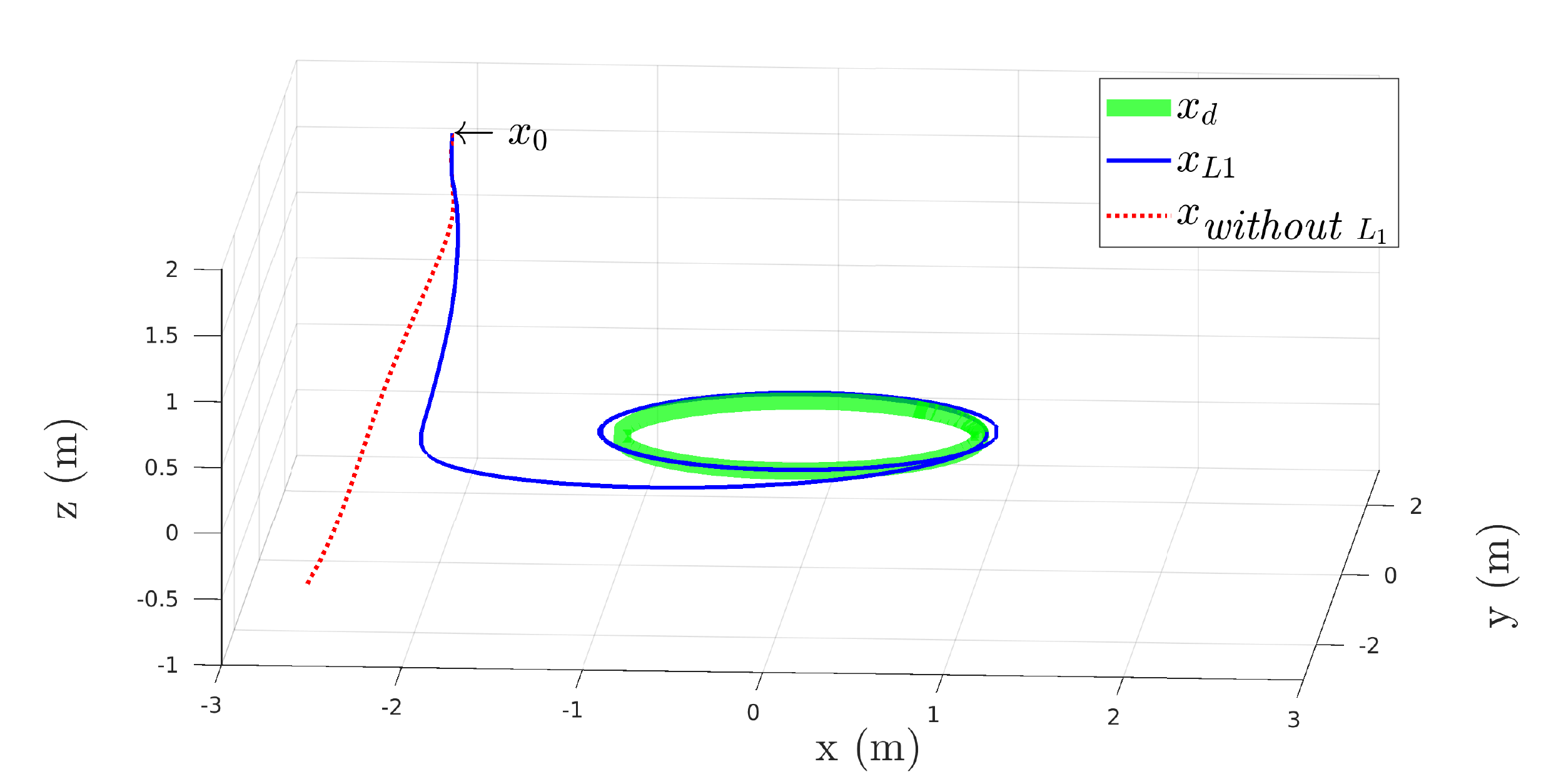}
\caption{\textit{Circular Trajectory - Case I}: Trajectory response for the two controllers defined earlier. Note that the geometric control \textit{without $L_1$} and with model uncertainty resulted in an unstable response (shown in red).}
\label{fig:circ_case_b_traj}
\end{figure}

\subsection{Geometric $L_1$ control vs Geometric control $\textit{without } L_1$}
\subsubsection{Case I: Constant External Disturbance $\theta_e$}
External disturbance of $\theta_t= [0.95,0.25,-0.5]^{T} \textit{Nm}$ in the body-frame is considered. $\theta_t$ along with an off-center added mass results in,
\begin{equation}
\label{eq:theta_e_sim}
\theta_e=[0.95,0.25,-0.5]^{T} + m_ag(r\times R^T\bm{e_3})\textit{Nm}.
\end{equation}
Note that, this value is unknown to the controller. The initial condition for the controller is almost an inverted case with a roll of $178^{\circ}$. Figure~\ref{fig:circ_case_b_errors} shows different errors in the system, including the errors between the \qmodel~and the \refsys~in the sub-figures \ref{fig_sim_1a} and \ref{fig_sim_1c}. Comparison between the controllers with $L_1$ and without $L_1$ are shown in the subfigures \ref{fig_sim_1b} and \ref{fig_sim_1d}. As shown in these figures, only $L_1$ control errors converge to zero.

Figure~\ref{fig:circ_case_b_theta} shows the uncertainty $\theta$ in the system, along with the predicted disturbance $\hat{\theta}$ and the filtered disturbance $C(s)\hat{\theta}$. As shown in the figure, the unfiltered uncertainty estimate is noisy while the filtered uncertainty estimate closely follows the actual uncertainty.

The quadrotor position, $x$, using the two different Controllers $L_1$, $\textit{without }L_1$ is shown in Figure~\ref{fig:circ_case_b_traj}. Note that the control $\textit{without }L_1$ fails to track the desired trajectory due to the large unknown model uncertainty in the system.

\subsubsection{Case II: Time varying disturbance $\theta_e(t)$}
Here, we consider a time varying disturbance along with model uncertainty. The following expressions denote the external disturbance used in this simulation example:
\begin{align}
\theta_t =& \frac{1}{2}(\cos(t)+0.5\cos(3t+0.23)\nonumber \\
&+0.5\cos(5t-0.4)+0.5\cos(7t+2.09), \label{eq:time_varying_theta}\\
\theta_e =&\theta_t + m_ag(r\times R^T\bm{e_3}).
\end{align}

\begin{figure}
    \centering
    \begin{subfigure}[t]{0.5\columnwidth}
        \centering
        \includegraphics[width=\columnwidth]{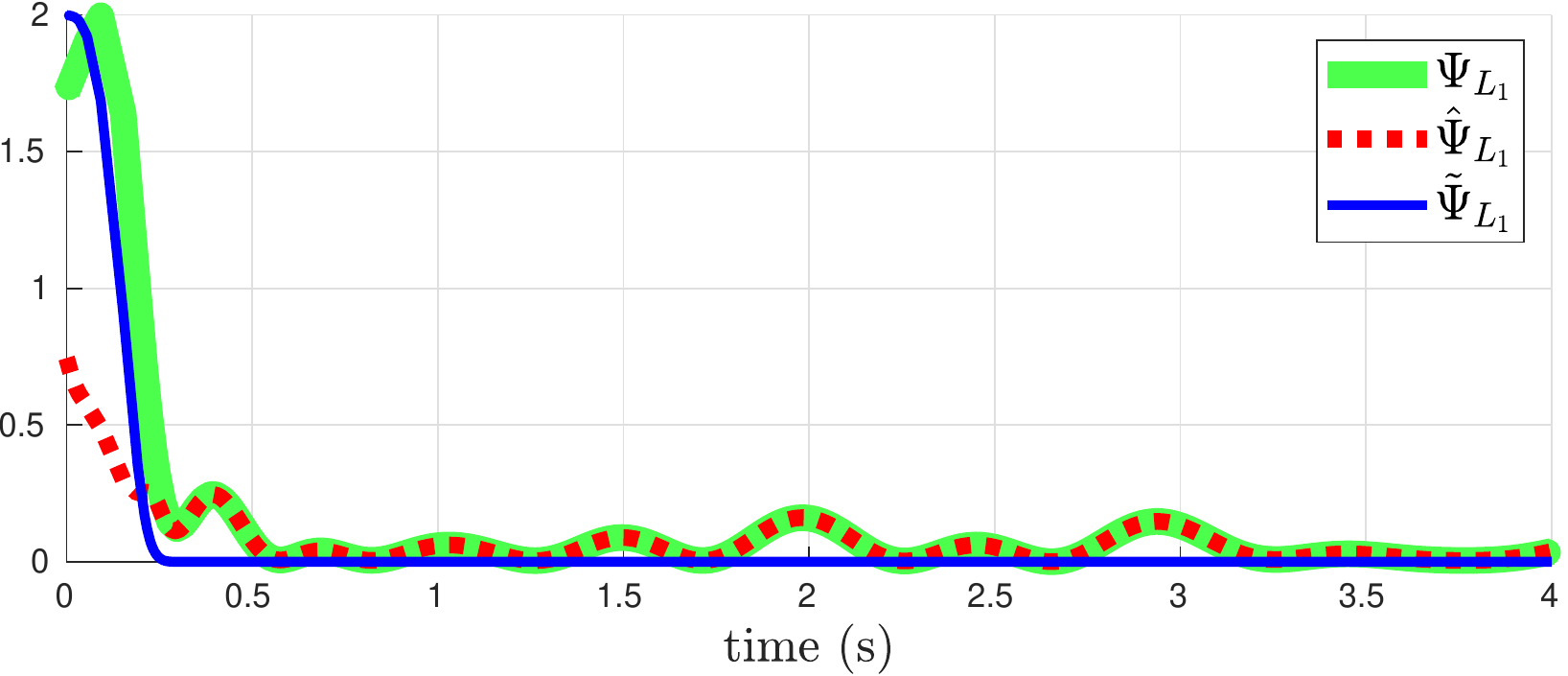}
        \caption{Configuration Errors for $L_1$: $\Psi\text{ vs }\hat{\Psi}\text{ vs }\tilde{\Psi}$}
        \label{fig_sim_2a}
    \end{subfigure}%
    ~ 
    \begin{subfigure}[t]{0.5\columnwidth}
        \centering
        \includegraphics[width=\columnwidth]{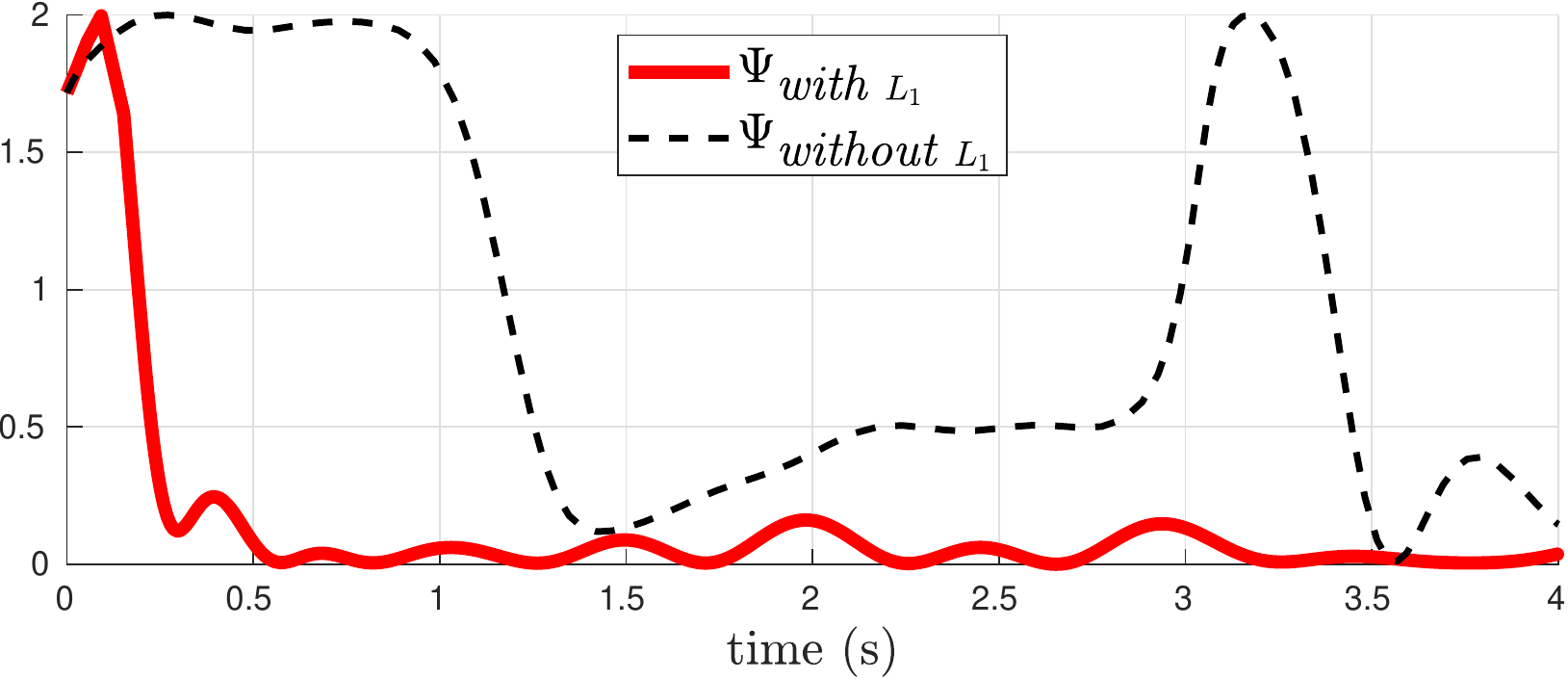}
        \caption{Configuration Errors: $L_1$ vs without $L_1$}
        \label{fig_sim_2b}
    \end{subfigure}
    
    \begin{subfigure}[t]{0.5\columnwidth}
        \centering
        \includegraphics[width=\columnwidth]{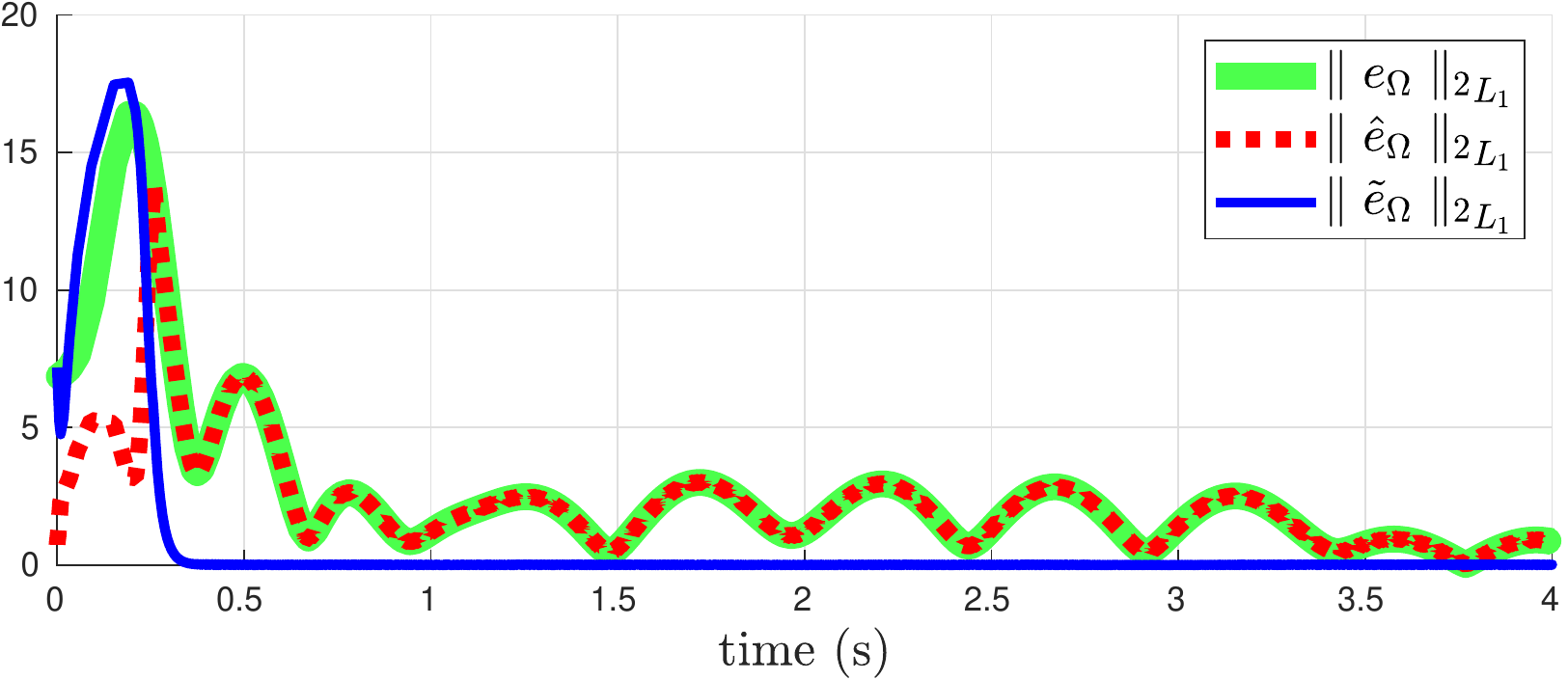}
        \caption{Angular Velocity Errors for $L_1$: $e_{\Omega}\text{ vs }\hat{e}_{\Omega}\text{ vs }\tilde{e}_{\Omega}$}
        \label{fig_sim_2c}
    \end{subfigure}%
    ~ 
    \begin{subfigure}[t]{0.5\columnwidth}
        \centering
        \includegraphics[width=\columnwidth]{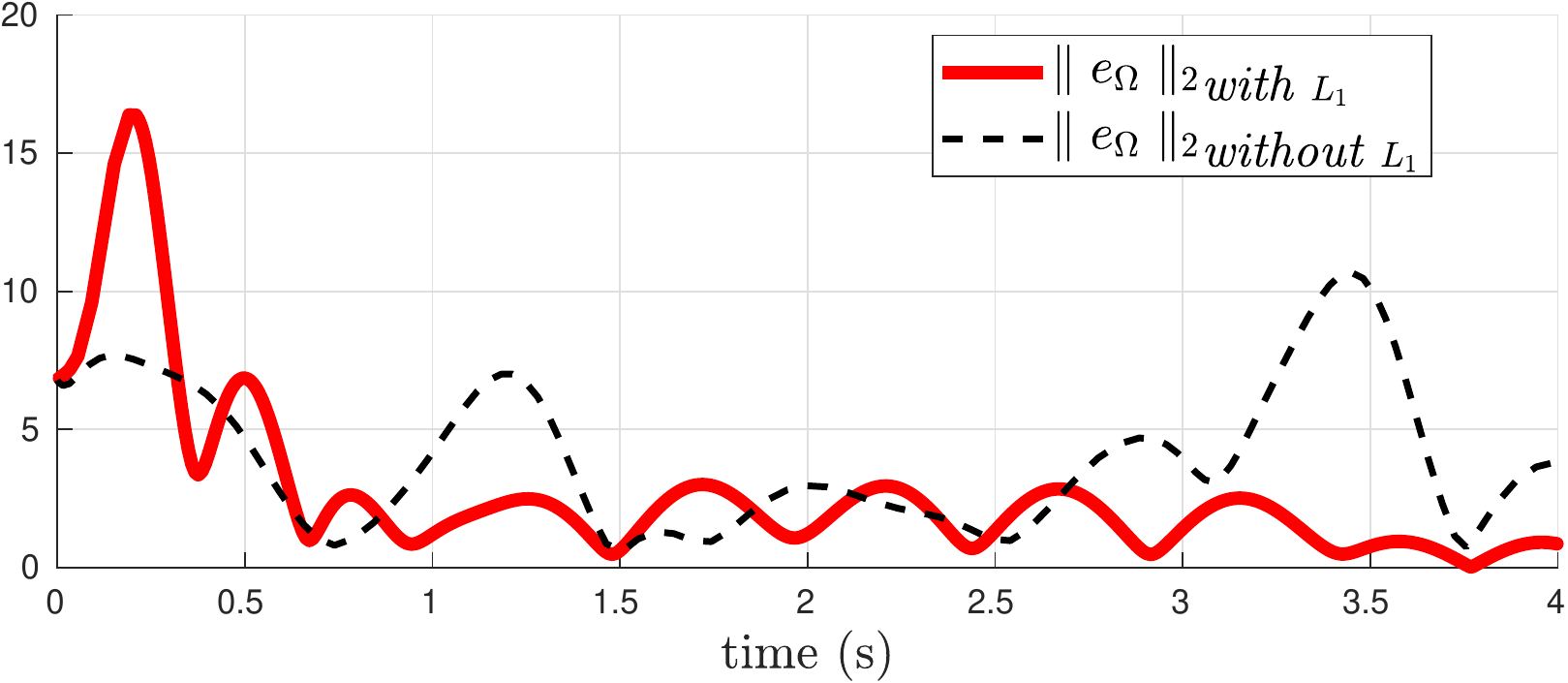}
        \caption{Angular Velocity Errors: $L_1$ vs without $L_1$}
        \label{fig_sim_2d}
    \end{subfigure}
    \caption{\textit{Circular Trajectory - Case II}: Comparison of errors between \qmodel\, and \refsys\, in the presence of time-varying external disturbances and model uncertainties (see \eqref{eq:time_varying_theta}, \eqref{eq:theta_e_1}). The \qmodel~errors and \refsys~errors did not converge to zero, however reached a bounded region about zero in case of geometric control \textit{with} $L_1$. Errors in case of geometric control \textit{without} $L_1$ are not-bounded as shown in \ref{fig_sim_2b} and \ref{fig_sim_2d}.}
    \label{fig:circ_case_c_errors}
\end{figure}
\begin{figure}
\centering
\includegraphics[width=1\columnwidth]{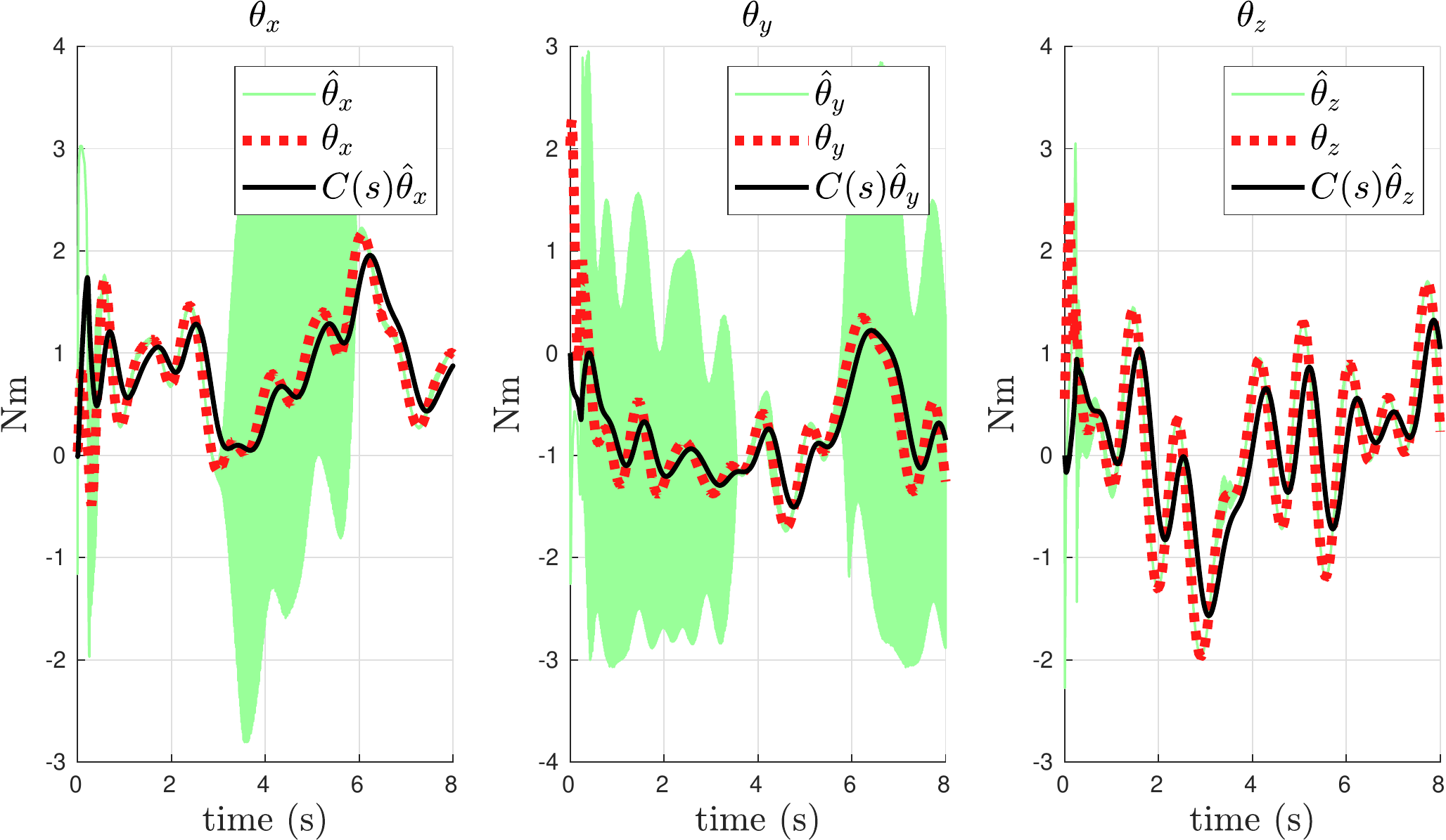}
\caption{\textit{Circular Trajectory - Case II}: Uncertainty in the system calculated using \eqref{eq:theta_def}. }
\label{fig:circ_case_c_theta}
\end{figure}
\begin{figure}
\centering
\includegraphics[width=1\columnwidth]{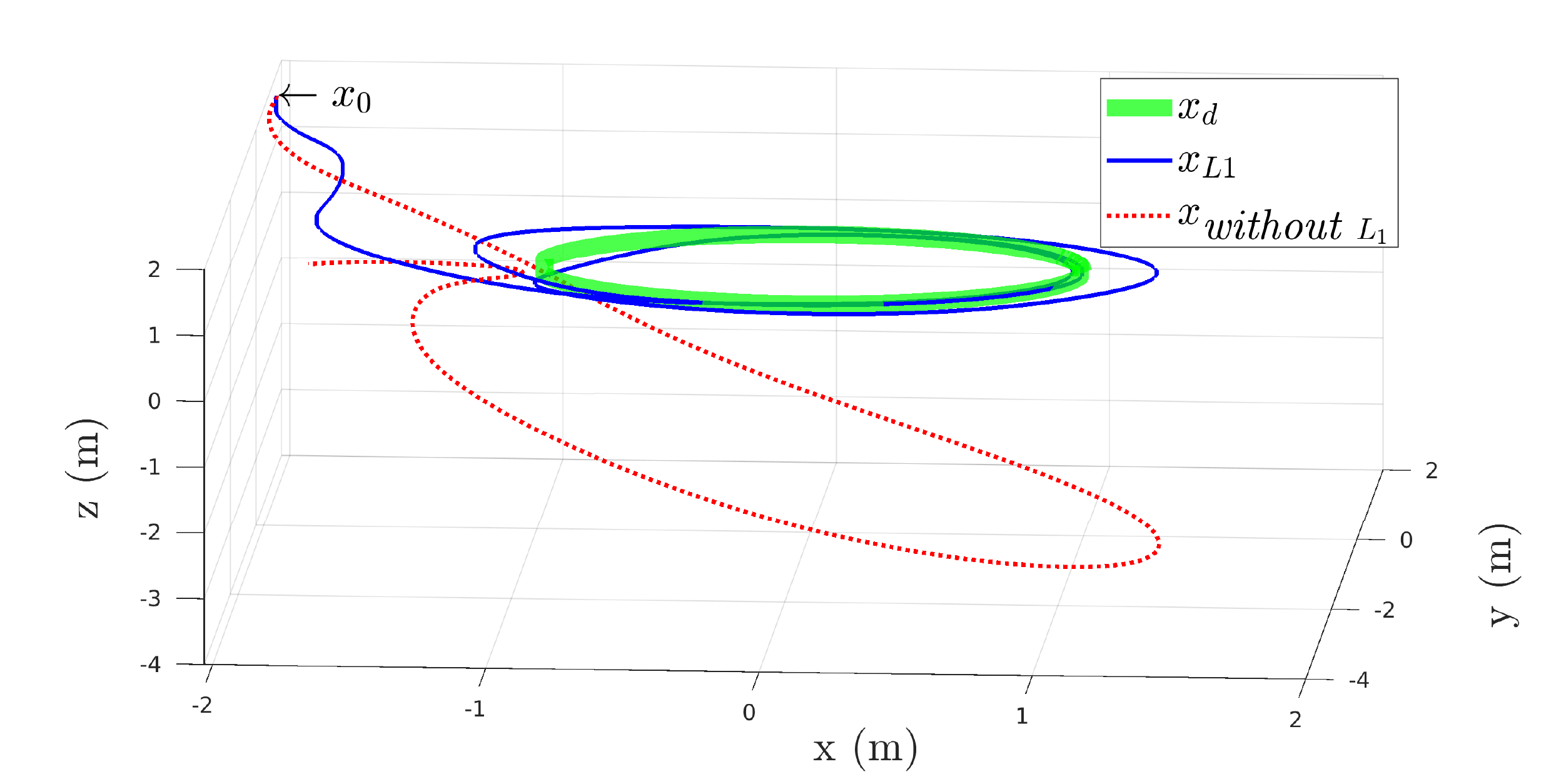}
\caption{\textit{Circular Trajectory - Case II}: Trajectory response for the two controllers defined earlier. Note that the geometric $L_1$ controller is able to track the desired trajectory even in the presence of time-varying disturbances (shown in blue).}
\label{fig:circ_case_c_traj}
\end{figure}
The uncertainty due to both the time-varying disturbance and model uncertainty and the corresponding estimation through $L_1$ adaptation is shown in Figure~\ref{fig:circ_case_c_theta}. Attitude tracking errors and resultant trajectory is shown in Figures~\ref{fig:circ_case_c_errors} and ~\ref{fig:circ_case_c_traj}. From sub-figures \ref{fig_sim_2a} and \ref{fig_sim_2c} it can be noticed that the errors between \qmodel\, and the \refsys\,,~$\tilde{\Psi}$ converge to zero very quickly, however the error between the \qmodel\, and the desired trajectory $\Psi$ doesn't completely converge to zero but reaches a bounded region about zero, this corresponds to the input-to-state stability discussed in Appendix C. As seen from Figure \ref{fig:circ_case_c_traj}, $L_1$ control performs better at tracking the trajectory even in the presence of time-varying disturbances. 
\highlight{\subsection{Geometric $L_1$ vs Euler $L_1$}}
\begin{figure*}
\centering
\includegraphics[width=2\columnwidth]{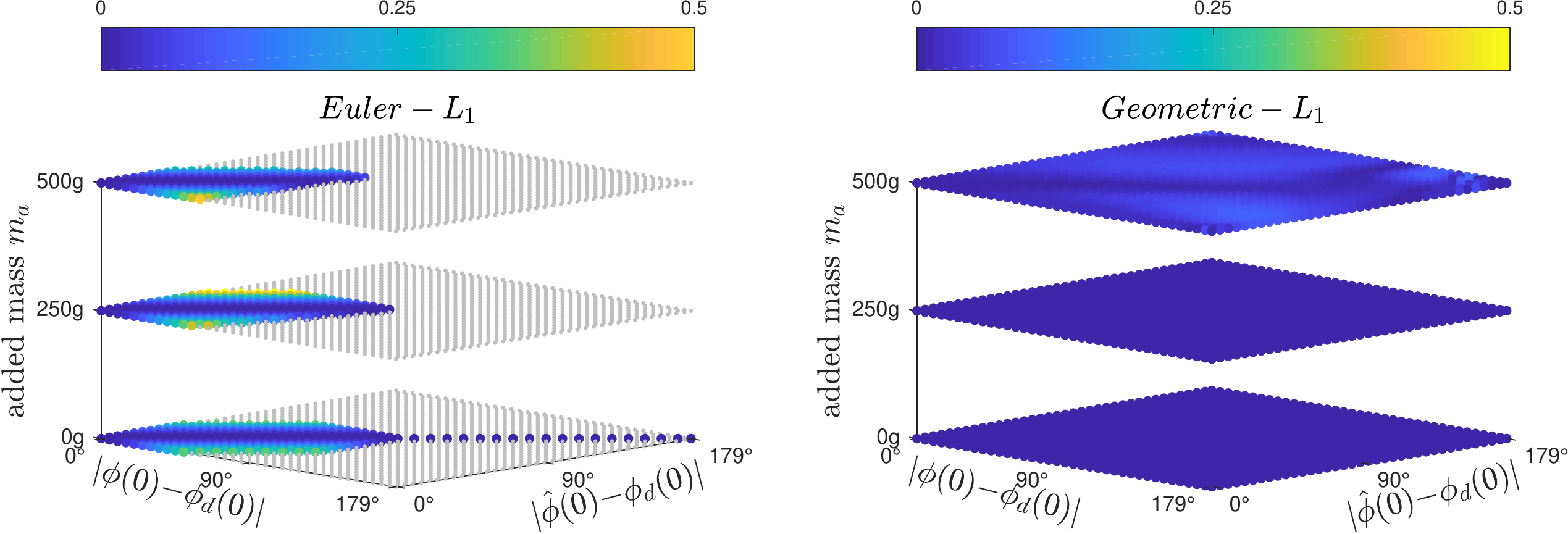}
\caption{\highlight{Comparison between Euler $L_1$ and Geometric $L_1$ \rebuttal{for various initialization errors}: Attitude configuration error, $\Psi$, after $3$ sec. of simulation for different added mass $\bm{m_a}$, different initial roll angle errors between reference model and the desired, $\hat \phi(0) - \phi_d(0)$, and between true model and desired, $\phi(0) - \phi_d(0)$. 
Initial roll angles for \qmodel~$\phi(0)$ and \refsys~$\hat{\phi}(0)$ are swept from $0^\circ$ to $179^\circ$ with $5^\circ$ increments, while the desired angles are $roll\equiv\phi_d\equiv~pitch{\equiv}~yaw{\equiv}0^\circ$. Note, that geometric $L_1$ works for all considered configuration errors while the Euler $L_1$ failed for larger configuration errors (shown in grey). With increased added mass, $\Psi$ increases and the range of errors for which the Euler $L_1$ control works decreases. } }%

\label{fig:evsg1}
\end{figure*}

\highlight{In order to establish the need for the geometric $L_1$ control, 
in this section we study and compare the initial condition response of Euler $L_1$ and geometric $L_1$. Similar to geometric $L_1$, an Euler $L_1$ control consists of a \refsys~and a \qmodel~, with the control moment computed using Euler angles and body-rates. The $L_1$ adaptation is implemented by considering a Lyapunov candidate for the Euler dynamics. 
In order to keep the comparison fair and avoid singularities, both controllers are simulated with the geometric rotational dynamics in \eqref{eq:R}-\eqref{eq:Omega}. }

\highlight{ 
We consider initial condition responses for several different initial roll angle errors between the \refsys~and the \textit{desired}, $\hat{\phi}(0)-\phi_d(0)$, and the \qmodel~and the \textit{desired}, $\phi(0)-\phi_d(0)$, for various added masses $\bm{m_a} \in \{0, 250g, 500g\}$.  We keep all other errors zero and consider zero desired angles, $roll\equiv\phi_d{\equiv}0^\circ,~pitch{\equiv}0^{\circ},~yaw{\equiv}0^\circ$. We ran the simulations varying $\hat{\phi}(0)-\phi_d(0)$ and $\phi(0)-\phi_d(0)$ from $0^\circ$ and $179^\circ$ in $5^\circ$ increments. The attitude configurations $\Psi$ after 3 seconds of simulations is shown in Figure~\ref{fig:evsg1} indicated by the colormap.  Each of the plots in Figure~\ref{fig:evsg1} illustrate $37 \times 37 \times 3 = 4107$ simulations.}

\highlight{We infer the following observations from the Figure,
\begin{enumerate}[{\it (i)}]
    \item Euler $L_1$ fails for larger attitude errors (shown in grey in the Figure), while geometric $L_1$ works for all considered configuration errors.
    \item Attitude configuration errors $\Psi$ increases with increasing model uncertainty (i.e., the added mass $\bm{m_a}$)
    \item The range of errors for which the Euler $L_1$ control works decreases with increased added mass $\bm{m_a}$
    \item Also, note that in the case of no added mass ($\bm{m_a}=0$), Euler $L_1$ works for all initial angles when \refsys~and \qmodel~ are initialized to the same angle, i.e., $\phi(0)=\hat{\phi}(0)$ corresponding to the diagonal line.
\end{enumerate}}

\begin{figure}
    \centering    
    \begin{subfigure}[t]{0.49\columnwidth}
        \centering
        \includegraphics[width=\columnwidth]{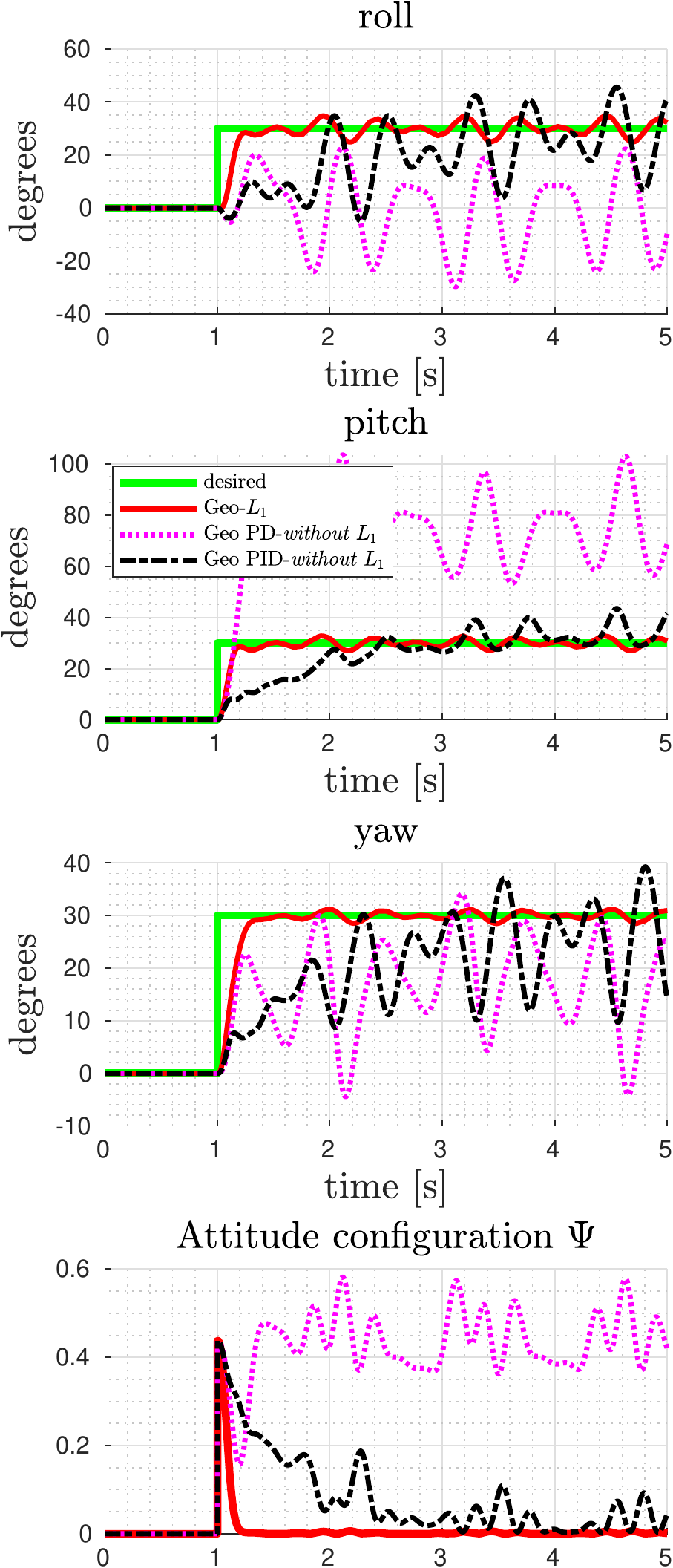}
        \caption{With Input Saturation}
        \label{fig_sim_2a}
    \end{subfigure}%
    ~ 
    \begin{subfigure}[t]{0.49\columnwidth}
        \centering
        \includegraphics[width=\columnwidth]{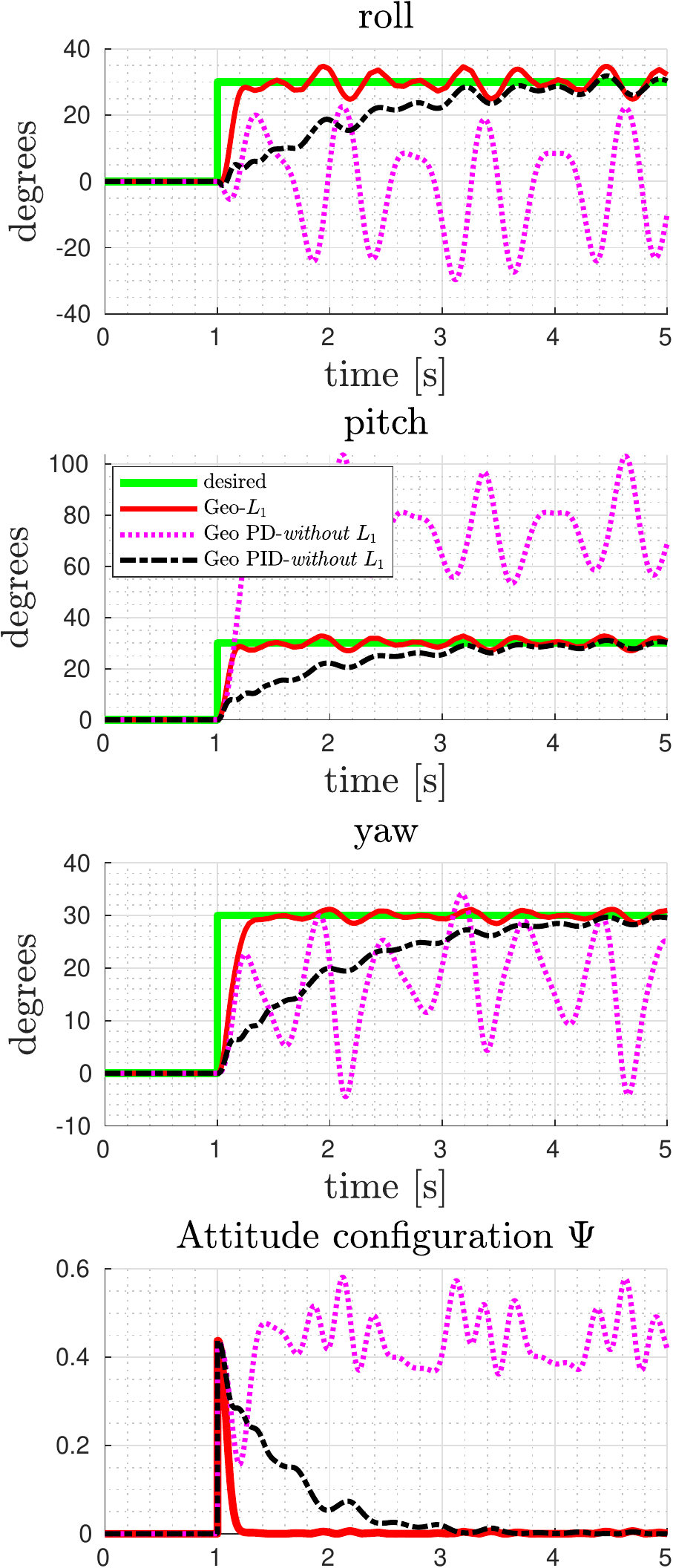}
        \caption{Without Input Saturation}
        \label{fig_sim_2b}
    \end{subfigure}
    \caption{\highlight{Step response comparison between different controls $(i)$ Geometric $L_1$, $(ii)$ Geometric PD $\textit{without }L_1$, $(iii)$ Geometric PID ($\textit{without }L_1$). Following desired Euler angles $roll{=}30^\circ$, $pitch{=}30^{\circ}$, $yaw{=}30^\circ$ are considered with added mass $\bm{m_a}=0.5kg$ and a time-varying disturbance. Note that PID control improves the steady-state performance when there is no input saturation, however, $L_1$ control has better transient performance in both cases.}} 
    \label{fig:step_responses}
\end{figure}

\highlight{\subsection{Step Input Response}
\label{sec:step-input-response}
In the section, we study the step response for geometric $L_1$ as well as other controllers.  We consider a step input change in desired angles $roll{=}30^\circ$, $pitch{=}30^{\circ}$, $yaw{=}30^\circ$, along with an uncertainty of $0.5kg$ added mass and a time-varying disturbance. We compare the performance of geometric $L_1$, geometric $PD$ \textit{without } $L_1$, and geometric $PID$ \cite{goodarzi2013geometric} (\textit{without } $L_1$). We compare these controllers both without and with an input saturation of $5$ Nm.
The PD gains for all three controllers are chosen the same.
The integral gain for the PID was increased to achieve the best performance for the case of no input saturation, and for the case of input saturation, the integral gain was chosen so that the input is just below the saturation in order to avoid integral windup.  The resulting step responses are shown in Figure~\ref{fig:step_responses}.
As seen from the Figure, the geometric $L_1$ control outperforms the geometric PD and also has a better transient performance than the geometric PID, irrespective of input saturation.  Furthermore, the geometric $L_1$ has a better steady-state performance compared to the geometric PID (with input saturation) and a similar steady-state performance compared to the geometric PID (without input saturation).

It is important to note that comparing different controllers is difficult as it involves tuning several parameters. For instance, PID control depends on the integral gain, while $L_1$ control depends on the adaptation gain and the cut-off frequency. In order to have a fair comparison, we applied the same input saturation to all controllers and studied their step response performance.}
	
	\section{Preliminary Experimental Results}
\label{sec:exp}
In this Section, we present the experimental results for the Geometric $L_1$ Adaptive Control developed in the previous sections.


\subsection{Setup}
\begin{figure}
\centering
\begin{subfigure}[t]{\linewidth}
	\includegraphics[width=1\columnwidth]{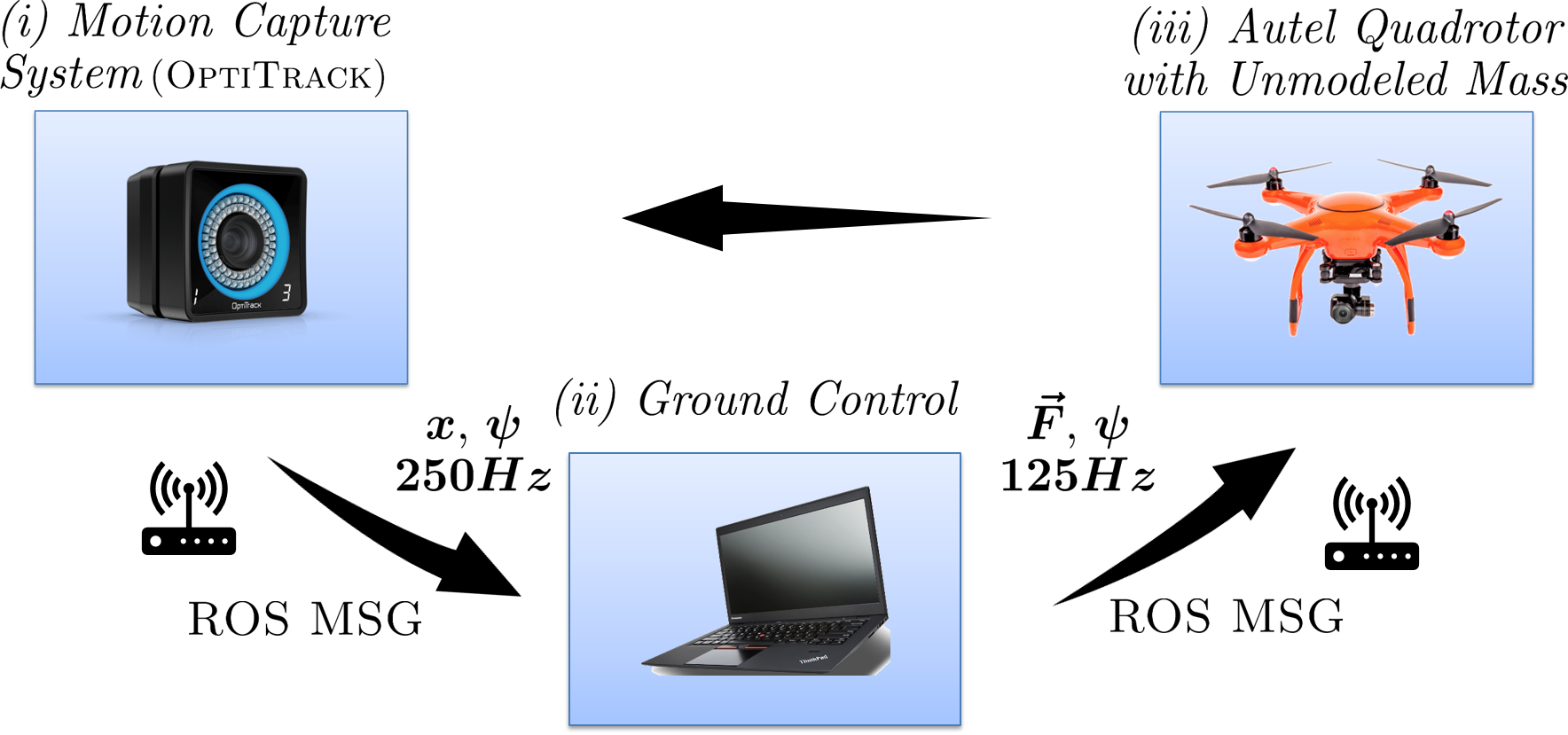}
	\caption{Setup used to validate the Geometric $L_1$ control consisting of  \textit{(i)} Motion capture system to track the pose of the quadrotor; \textit{(ii)} Ground control with the position control for the quadrotor; and \textit{(iii)} Autel quadrotor with on-board Geometric Attitude control with $L_1$ adaptation.}
	\label{fig:exp-setup}
\end{subfigure}
\begin{subfigure}[t]{\linewidth}
    \includegraphics[width=\linewidth]{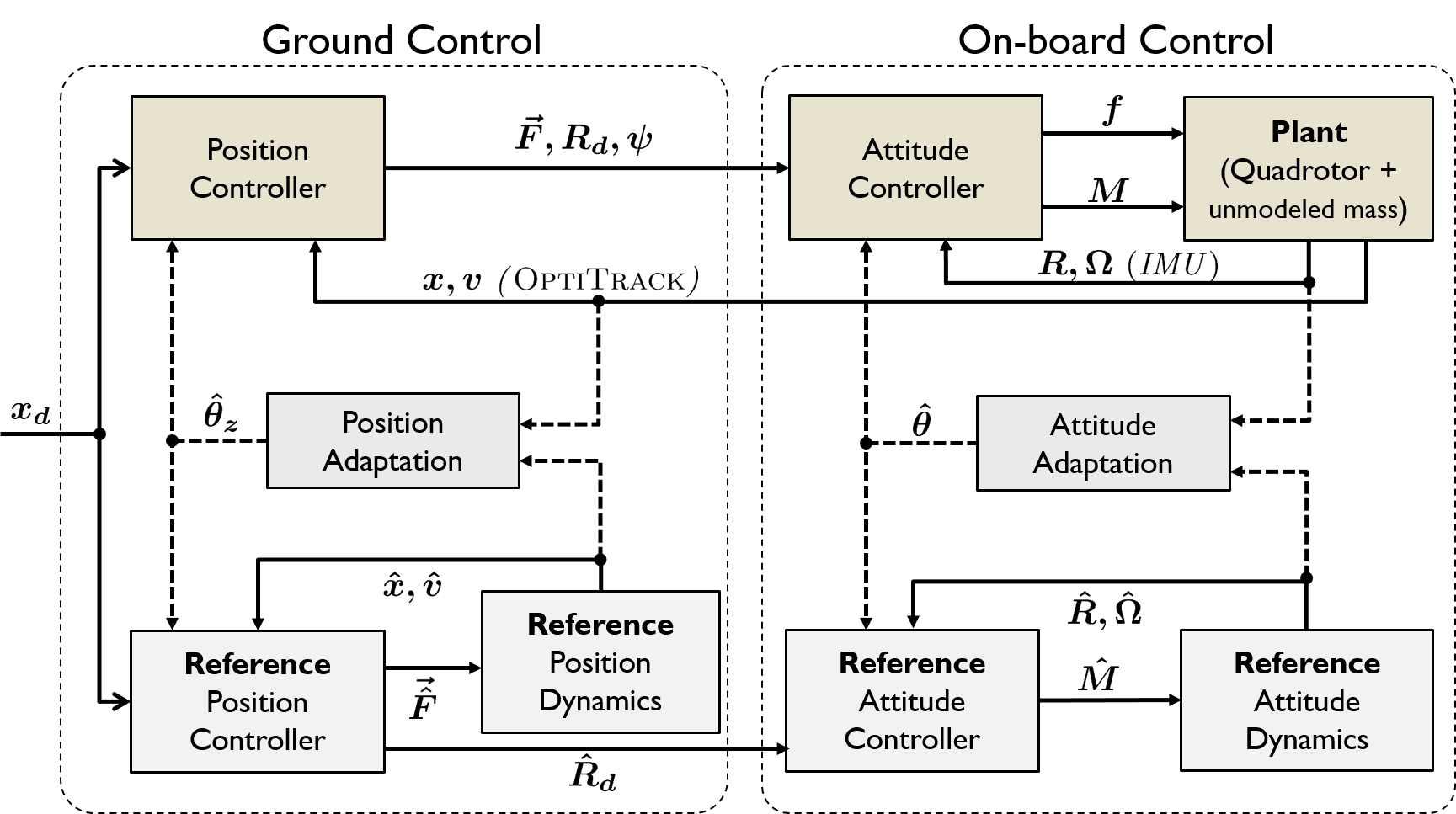}
    \caption{Control architecture used in the experiments showing the ground and on-board controls.}
    \label{fig:ctrl_struct}
\end{subfigure}
\caption{Experimental setup and the control architecture to track the quadrotor trajectory with $L_1$ adaptation to estimate the uncertainty $\theta$.}
\label{fig:exp-setup_full}
\end{figure}

The experiments are conducted using the Autel-X star quadrotor equipped with a Raspberry Pi 3 based Navio-2. A ROS node on Raspberry Pi 3 runs the on-board attitude control at rate of $1kHz$. Figure~\ref{fig:exp-setup} illustrates the experimental setup used in this paper. A motion capture system Optitrack is used to estimate the pose, velocity and yaw of the quadrotor at $250Hz$. The Inertial Measurement Unit (IMU) on the Navio2 is used to estimate the body-attitude and body-rates. A Lenovo-Thinkpad with Ubuntu 14.04 and ROS constitutes the ground control. A ROS node runs the position control on the ground control and communicates with the onboard control through WiFi at $125Hz$.

\highlight{To generate disturbances and model uncertainty to the quadrotor system, a weight of $\bm{m_a}= 0.2kg$ is rigidly attached to the fuselage of the quadrotor at approximately $[0.1, 0.1, 0.1]^Tm$ in the body-frame. This additional mass is unknown to the controller(s) and constitutes both model uncertainty, with true inertial ${\Jtrue}$ (see \eqref{eq:true_inertial_calc}) and disturbance in the system; added mass generates a moment due to gravity about the quadrotor center-of-mass, given as,
    \begin{align*}
        \theta_m = \bm{m_a}g(r\times R^T\bm{e_3}). 
    \end{align*}
Since, this moment varies with $R$, especially for a circular trajectory it would result in a time-varying disturbance.   }

\subsection{Control Architecture}
The control technique used to track the desired trajectory using $L_1$ Geometric Attitude controller consists of $(i)$ a position controller, along with a reference model for position dynamics and $(ii)$ an attitude controller with a reference attitude dynamics.
Figure~\ref{fig:ctrl_struct} illustrates the control architecture used in the experiments.
From the desired trajectory, positions and its higher derivatives are calculated and used to calculate the desired states and the feed forward inputs using differential flatness~\cite{mellinger2011minimum}. From the desired position, velocity and acceleration, the desired thrust vector $\vec{F}$ is calculated. A reference position dynamics is simulated using Euler integration and the uncertainty in the position dynamics is estimated and is compensated in the thrust calculation. Based on this thrust vector information, the desired attitude is calculated on-board and is used to compute the desired moments for the quadrotor $M$ and reference model $\hat{M}$.
$\hat{M}$  is used in the reference attitude dynamics simulation, achieved through on-board Euler integration. The scalar thrust, $f$, is calculated using the thrust vector, $\vec{F}$, and the attitude $R$ using \eqref{eq:f}. The moment, $M$, along with the thrust, $f$, is used to generate the desired angular speed for the motors as presented in \cite{mahony2012multirotor}.

\subsection{Results}

To validate the developed controller, we show the tracking performance of the quadrotor (with the weight $\bm{m_a}$ attached) with geometric $L_1$ adaptive control and geometric PD control in \eqref{eq:Mactual} \highlight{without $L_1$ adaptation (referred to as \textit{without} $L_1$)}.
\highlight{
\begin{remark}
Note that the position controller in the ground control (see Figure~\ref{fig:ctrl_struct}) is same for both $L_1$ and \textit{without} $L_1$ controllers (i.e, we are comparing only the attitude control).
\end{remark} }

The different experimental system parameters for the quadrotor are,
\begin{gather*}
\bm{m} = 1.129kg,~\bm{m_a} = 0.2kg,\,\highlight{\bm{r} = [
0.1, 0.1, 0.1]^Tm}
\end{gather*}
{\scriptsize \begin{gather*}
\bm{J} = \begin{bmatrix}
6.968 & -0.02909 & -0.2456 \\ -0.02909 & 6.211 & 0.3871 \\ -0.2456 & 0.871 & 10.34
\end{bmatrix}10^{-3}kgm^2.
\end{gather*}}
We show the performance for Hover and a Circular trajectory tracking in the following subsections. \highlight{Video showcasing the experiments can  can be found at this link, \href{https://youtu.be/nBDDxpkz6Pg}{https://youtu.be/nBDDxpkz6Pg}. }

\subsubsection{\bf Hover}
\begin{figure}
\begin{subfigure}[t]{\linewidth}
    \includegraphics[width=\linewidth]{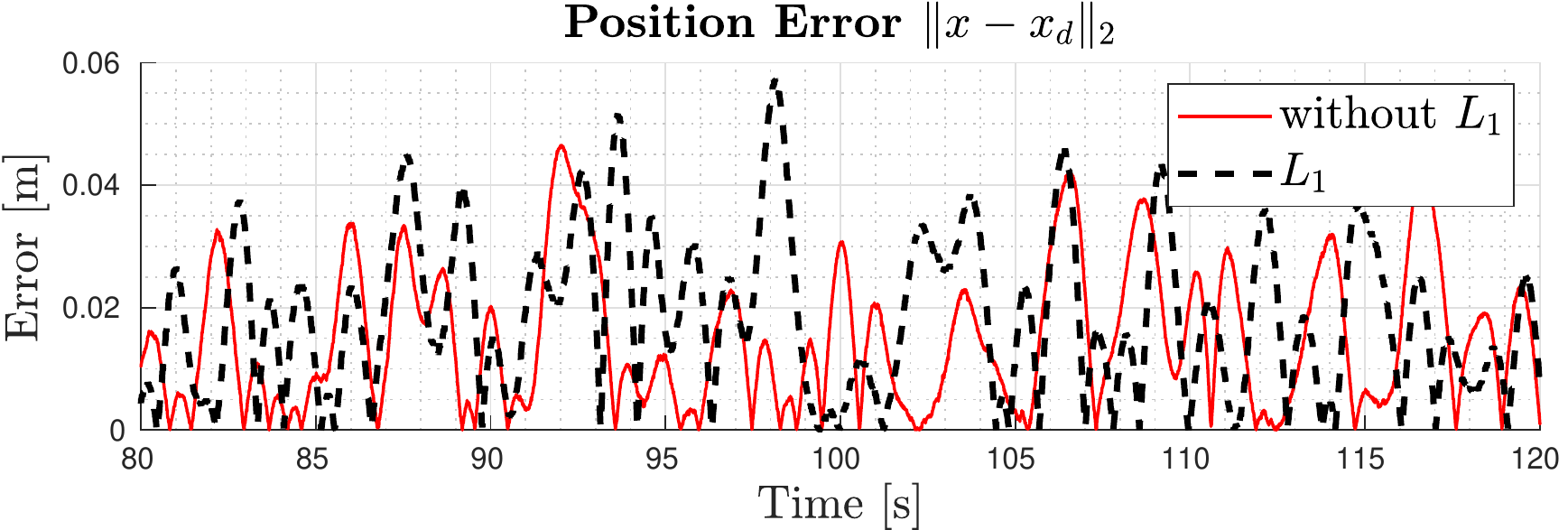}
    \caption{Position error $\|x-x_d\|_2$(m)}
    \label{fig:case2pos}
\end{subfigure}

\begin{subfigure}[t]{\linewidth}
    \includegraphics[width=\linewidth]{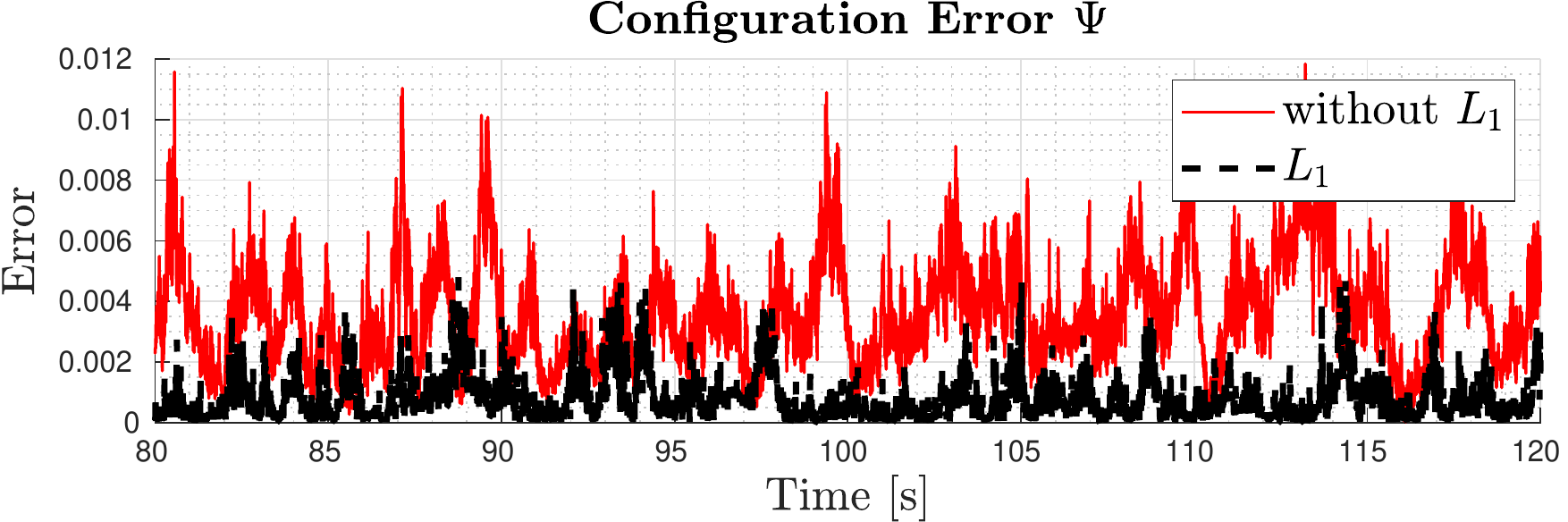}
    \caption{Attitude configuration error $\Psi$}
    \label{fig:case2Psi}
\end{subfigure} 
\caption{\textit{Hover}: Tracking performance with and without $L_1$ controller on quadrotor with attached weight. (a) shows the position error, and (b) the attitude error.}
\label{fig:hover_exp}
\end{figure}

The desired states for the Hover are, \highlight{$x_d=[0, -1, 1]^T$, $v_d = [0, 0, 0]^T$, $R=I_{3\times 3}$ and $\Omega_d=[0, 0, 0]^T$}. Figure~\ref{fig:hover_exp} presents the tracking performance in terms of (a) position error, and (b) the attitude configuration error for Hover. The mean and standard-deviation for these errors are presented in Table~\ref{table:3_errors}.

\subsubsection{\bf Circular Trajectory}
\begin{figure}
\begin{subfigure}[t]{\linewidth}
    \includegraphics[width=\linewidth]{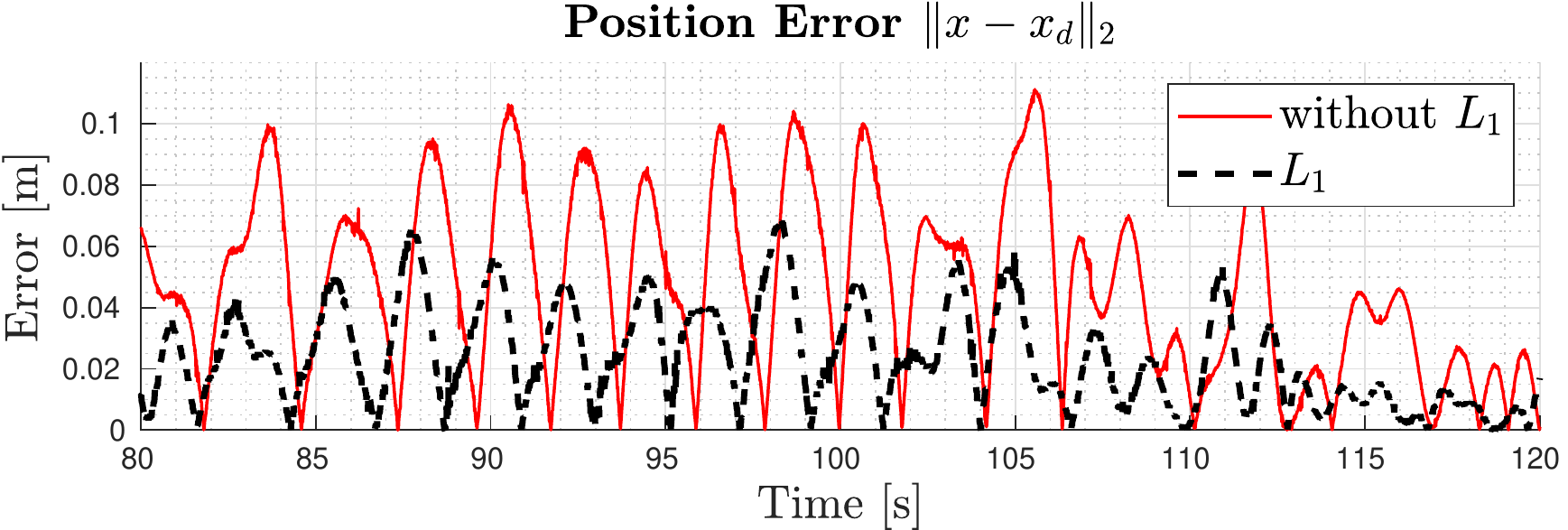}
    \caption{Position error $\|x-x_d\|_2$(m)}
    \label{fig:case2pos}
\end{subfigure}

\begin{subfigure}[t]{\linewidth}
    \includegraphics[width=\linewidth]{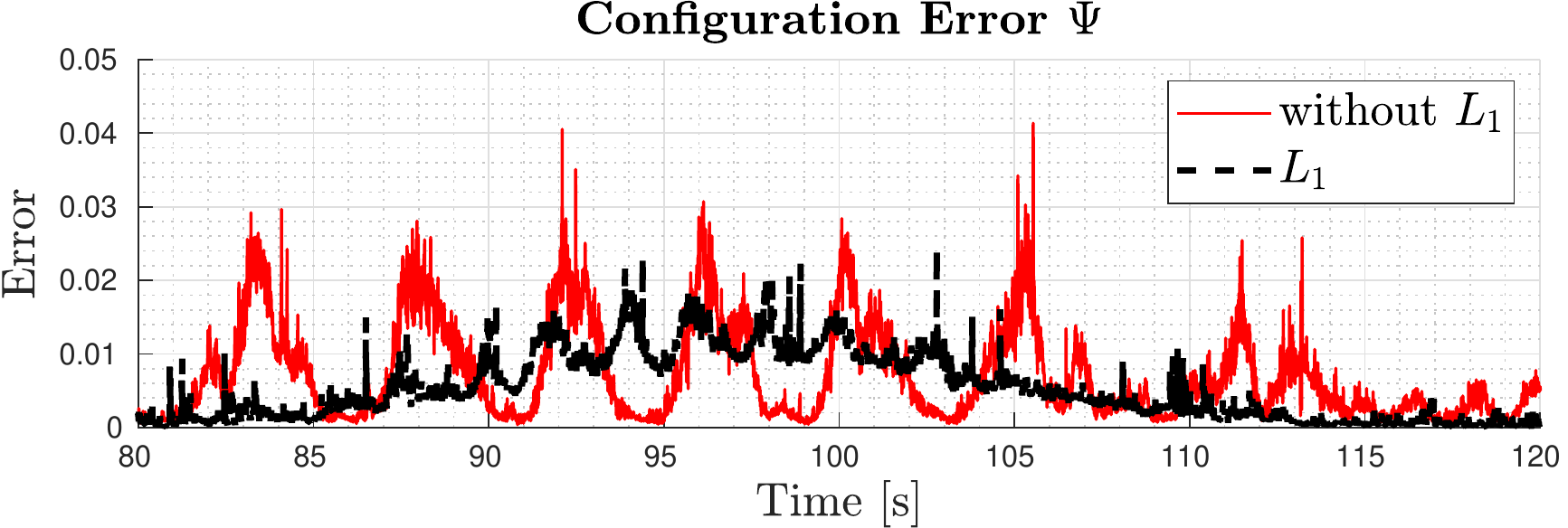}
    \caption{Attitude configuration error $\Psi$}
    \label{fig:case2Psi}
\end{subfigure} 
\caption{\textit{Circular Trajectory}: Tracking performance with and without $L_1$ controller on quadrotor with attached weight. (a) shows the position error, (b) the attitude error.}
\label{fig:circle_exp}
\end{figure}

We chose the following circular trajectory in flat-outputs:

\begin{gather}
\label{eq:trajectory_flat}
\highlight{x = [ x_0{+}\rho\cos(\omega (t)), y_0{+}\rho\sin(\omega (t)),z_0]^T, \psi(t) = 0} 
\end{gather}
where, $\omega(t) = \frac{2\pi a}{b+\exp^{-c(t-t_0)}}$ and $[x_0,\,y_0,\,z_0]$ is the center of the circle and \highlight{$\rho$} is the radius of the circle. The parameters used in the experiment are \highlight{$\rho=1$, $x_0=y_0=0,~z_0=1.5$} $a=13$, $b=1$, $c=0.1$ and $t_0 = 80$. These values result in a circular trajectory with increasing speed to a maximum of $2m/s$ and then decreasing to zero. Tracking performance for this trajectory is presented in Figure~\ref{fig:circle_exp} and the mean and standard-deviation of the errors are presented in Table~\ref{table:3_errors}.

\begin{table}
\scriptsize
\centering
\begin{tabular}{| l | c |  c | c|} 
\hline 
\multirow{2}*{Experiment} & \multirow{2}*{Error} & Without $\bm{L_1}$ & With $\bm{L_1}$\\
						& 						& mean [std-dev] & mean [std-dev]\\	
\hline \hline 
\multirow{2}*{Hover} & $\| x-x_d\|_2 $ & $0.0152$ [$0.0118$] & $0.0182$ [$0.0132$]\\
					& $\Psi $ & $0.0055$ [$0.0462$] & $0.0021$ [$0.0452$]\\
\hline
\multirow{2}*{\parbox{1.5cm}{Circular \\ Trajectory}} & $\| x-x_d\|_2$ & $0.0250$ [$0.0220$]  & $0.0126$ [$0.0120$] \\
					& $\Psi $ & $0.0056$ [$0.0428$]  & $0.0027$ [$0.0423$] \\
\hline
\end{tabular}
\caption{Mean and standard-deviation for position and attitude errors for Hover and Circular Trajectory with and without $L_1$ adaptation.}
\label{table:3_errors}
\end{table}

\begin{remark}
As it can be noticed from Figures~\ref{fig:hover_exp} \& ~\ref{fig:circle_exp} and Table~\ref{table:3_errors}, controllers with and without $L_1$ have similar performance in case of Hover, while the $L_1$ controller shows approximately a factor of two better performance in the circular trajectory case. This could be attributed to the fact that the states and especially the body-angular velocity has an effect on the uncertainty as can be seen from \eqref{eq:Om_real_3} \& \eqref{eq:theta_def}. In particular, in case of hover the angular velocity $\Omega$ is close to zero resulting in a small uncertainties $\theta_m \approx 0$ while in case of the circular trajectory, the feed-forward angular velocity results in significant uncertainty $\theta$. 
\end{remark}

	\section{Conclusion}
\label{section: con}
We address the problem of disturbances and model uncertainties in cases of quadrotor controlled using geometric control techniques, where coordinate free dynamics are used to avoid singularities. We develop the error dynamics with uncertainties in the system and a reference model without any uncertainty. Attitude tracking errors are defined between the robot model and the reference model. We develop control for the robot model and the reference model with $L_1$ adaptation to estimate the uncertainty in the robot model. Control Lyapunov candidate is defined to show that the proposed control strategy results in an exponential stability between the robot model and reference model, beyond a small bound about the origin. The bounded region is inversely  proportional to the chosen adaptation gain $\Gamma$. Numerical simulations are presented to validate the control in presence of disturbances and model uncertainties.
Experimental results are presented to show the performance of the $L_1$ adaptive control.

\appendix
\label{appendix:1}

\subsection*{\highlight{Appendix A: Useful Mathematical Background}}
\label{sec:appendixA}
The Cross-map and Vee-map were defined in Section~\ref{section: quad_dynamics}. Following are few properties that are useful in the rest of the appendix. For any $x \in \mathbb{R}^{3}$, $A \in \mathbb{R}^{3 \times 3}$ and $R \in SO(3)$, we have,

\begin{align}
tr[x^{\times}A]&=tr[Ax^{\times}] \label{eq:hatprop1}\\
& =\frac{1}{2}tr[x^{\times}(A-A^{T})] \label{eq:hatprop1b} \\
&=-x^{T}(A-A^{T})^{\vee}, \label{eq:hatprop1c}\\
Rx^{\times}R^{T}& =(Rx)^{\times}, \label{eq:hatprop2}\\
\hatmap{x}A + A^T\hatmap{x} &= \hatmap{(\{tr[A]I-A\}x)}. \label{eq:hatprop3}
\end{align}

\subsection*{\highlight{Appendix B: Attitude Error Dynamics}}
As described in earlier sections, attitude is represented using a rotation matrix $R \in SO(3)$ and body-angular velocity $\Omega\in T_RSO(3)$. Errors between different attitude and angular velocities namely, true system, reference system and desired trajectory are different and defined in \eqref{eq:eRdef} \& \eqref{eq:eOmdef}, \eqref{eq:eRhat} \& \eqref{eq:eOmhat} and \eqref{eq:eRtilde} \& \eqref{eq:eOmtilde}. We compute the derivatives of the errors between \qmodel~ and \refsys, $(\tilde{e}_{R}$ and $\tilde{e}_{\Omega})$ as shown below.

\begin{align}
\label{eq:eRtdot}
\dot{\tilde{e}}_{R} &=\frac{1}{2}\Big(\frac{d}{dt}(R^T\hat{R})-\frac{d}{dt}(\hat{R}^TR)\Big)^{\vee},
\end{align}
\begin{remark}
We have,
\begin{align}
\frac{d}{dt}(R^T\hat{R})&= \dot{R}^T\hat{R} + R^T\dot{\hat{R}} \nonumber \\
&=(R\hatmap{\Omega})^T\hat{R} + R^T(\hat{R}\hatmap{\hat{\Omega}})\nonumber \\
&=R^T\hat{R}(-\hat{R}^TR\hatmap{\Omega}R^T\hat{R} + \hatmap{\hat{\Om}})
\end{align}
Using \eqref{eq:hatprop2}, we get
\begin{align}
\frac{d}{dt}(R^T\hat{R})=R^T\hat{R}\hatmap{(\hat{\Om} - \hat{R}^TR\Om)}=R^T\hat{R}\hatmap{\tilde{e}}_{\Omega},
\end{align}
and consequently,
\begin{align}
\frac{d}{dt}(\hat{R}^TR)&=-\hatmap{\tilde{e}}_{\Omega}\hat{R}^TR.
\end{align}
\end{remark}
Therefore,
\textbf{\begin{equation}
\dot{\tilde{e}}_{R} = \frac{1}{2}( R^{T}\hat{R}\tilde{e}_{\Omega}^{\times} + \tilde{e}_{\Omega}^{\times}\hat{R}^{T}R)^{\vee},
\end{equation}}
Using the property in \eqref{eq:hatprop3} we get,
\begin{equation}
\label{eq:eRtdot}
\dot{\tilde{e}}_{R} = \frac{1}{2}(Tr[\hat{R}^{T}R]I-\hat{R}^{T}R)\tilde{e}_{\Omega} \eqqcolon C(R^{T}\hat{R})\tilde{e}_{\Omega},
\end{equation}
where it is shown in \cite{lee2011geo} that the function $C(R^T\hat{R})$ satisfies the property $\|C(R^{T}\hat{R})\tilde{e}_{\Omega}\|_{2} \leq 1$ for any rotation matrix in SO(3).

Next, from \eqref{eq:eOmtilde} we have,
\begin{align}
\tilde{e}_{\Om}&=\hat{\Omega}-\hat{R}^TR\Omega\\
&=\hat{\Omega}-\hat{R}^TR_d\Omega_d + \hat{R}^TR_d\Omega_d-\hat{R}^TR\Omega\nonumber \\
&=(\hat{\Omega}-\hat{R}^TR_d\Omega_d)-\hat{R}^TR(\Omega-{R}^TR_d\Omega_d)\nonumber\\
\implies \tilde{e}_{\Omega}&=\hat{e}_{\Omega}-\hat{R}^TRe_{\Omega}.
\end{align}

The derivative of $\tilde{e}_{\Omega}$ then is,
\begin{align}
\dot{\tilde{e}}_{\Omega} & = \dot{\hat{e}}_{\Om}-\big( -\hatmap{\tilde{e}}_{\Om}\hat{R}^TRe_{\Om} + \hat{R}^TR\dot{e}_{\Om}\big), \label{eq:dot_e_Omega}\\
\implies J\dot{\tilde{e}}_{\Omega} & = J\dot{\hat{e}}_{\Om}-J\big( -\hatmap{\tilde{e}}_{\Om}\hat{R}^TRe_{\Om} + \hat{R}^TR\dot{e}_{\Om}\big). \label{eq:Jdot_e_Omega}
\end{align}
From \eqref{eq:eOm_real_final} and \eqref{eq:deOmhat}, we get
\begin{align}
J\dot{\tilde{e}}_{\Omega} & = \hat{\mu} -J\big( 
{-\hatmap{\tilde{e}}_{\Om}\hat{R}^TRe_{\Om}} + \hat{R}^TRJ^{-1}\big[\mu + \theta\big]\big).
\end{align}

From the control moments defined in \eqref{eq:controlMomentM_first} - \eqref{eq:controlMomentM} we have,
\begin{align}
J\dot{\tilde{e}}_{\Omega}  =& \hat{\mu} -J\Big( 
{-\hatmap{\tilde{e}}_{\Om}\hat{R}^TRe_{\Om}} \nonumber\\
&+\hat{R}^TRJ^{-1}\big[JR^T\hat{R}J^{-1}
(\hat{\mu}+\tilde{k}_R\tilde{e}_R+\tilde{k}_{\Omega}\tilde{e}_{\Omega}) \nonumber \\
&-R^T\hat{R}\hat{\theta} + JR^T\hat{R}\hatmap{\tilde{e}}_{\Omega}\hat{R}^TRe_{\Omega}+ \theta\big]\Big)
\end{align}

\begin{align}
& = -\tilde{k}_R\tilde{e}_R - \tilde{k}_{\Om}\tilde{e}_{\Om} +  J\hat{R}^TRJ^{-1}R^T\hat{R}(\hat{\theta} -\hat{R}^TR\theta).
\end{align}
Thus, we finally have,
\begin{align}
J\dot{\T{e}}_{\Omega}&=-\tilde{k}_R\tilde{e}_R - \tilde{k}_{\Omega}\tilde{e}_{\Omega} + \Rho\tilde{\theta} \label{eq:final_deOmtilde}
\end{align}
where $\Rho$ and $\tilde{\theta}$ are defined in \eqref{eq:Rho_def}, \eqref{eq:thetatilde} respectively.

The time-derivative of the configuration error function $\Tilde{\Psi}$ (see \eqref{eq:Psitilde}) can be computed and simplified using \eqref{eq:hatprop2} as follows,
\begin{align}
\dot{\tilde{\Psi}} &= -\frac{1}{2}Tr[R^{T}\hat{R}(\hat{\Omega}^{\times} - \hat{R}^{T}R\Omega^{\times}R^{T}\hat{R})]\\
&= -\frac{1}{2}Tr[R^{T}\hat{R}(\hat{\Omega} - \hat{R}^{T}R\Omega)^{\times}] \\
&= -\frac{1}{2}Tr[R^{T}\hat{R}\tilde{e}_{\Omega}^{\times}].
\end{align}
Then using \eqref{eq:hatprop1c} we get,
\begin{equation}
\label{eq:Psitdot}
\dot{\tilde{\Psi}} = \frac{1}{2} \tilde{e}_{\Omega}^{T}(R^{T}\hat{R}-\hat{R}^{T}R)^{\vee} = \tilde{e}_{\Omega}^{T}\tilde{e}_{R} = \tilde{e}_{\Omega} \cdot \tilde{e}_{R}.
\end{equation}

\subsection*{\highlight{Appendix C: Lyapunov Function Candidate}}
\label{appendix:LFC}
\textbf{Proof for Theorem \ref{prop:1}}: To show the \textit{exponential input-to-state stability} of the attitude errors, $(\tilde{e}_R,\tilde{e}_{\Omega})$, we will consider the following control Lyapunov candidate function,
\begin{equation}
\label{eq:V}
V=\frac{1}{2}\tilde{e}_{\Omega} \cdot J\tilde{e}_{\Omega}+\tilde{k}_{R}\tilde{\Psi}(\hat{R},R)+c\tilde{e}_{R} \cdot \tilde{e}_{\Omega}+\frac{1}{2}\tilde{\theta}^{T}\Gamma^{-1}\tilde{\theta},
\end{equation}
where the adaptation gain, $\Gamma$, is a symmetric, positive-definite matrix and $\tilde{k}_R$ is a positive number.
Taking the derivative of $V$ gives,
\begin{align}
\dot{V} =& \tilde{e}_{\Omega} \cdot J\dot{\tilde{e}}_{\Omega} +\tilde{k}_{R}\dot{\tilde{\Psi}}(\hat{R},R)\nonumber \\
&+c\dot{\tilde{e}}_{R} \cdot \tilde{e}_{\Omega} +c\tilde{e}_{R} \cdot \dot{\tilde{e}}_{\Omega} + \tilde{\theta}^{T}\Gamma^{-1}\dot{\tilde{\theta}}. \label{eq:Vdot}
\end{align}

Substituting equations \eqref{eq:eRtdot}, \eqref{eq:final_deOmtilde}, and \eqref{eq:Psitdot} for $\dot{\tilde{e}}_{R}$, $\dot{\tilde{e}}_{\Omega}$ and $\dot{\tilde{\Psi}}$ respectively, we get,
\begin{align}
\dot{V} =& \tilde{e}_{\Omega} \cdot(-k_{R}\tilde{e}_{R}-k_{\Omega}\tilde{e}_{\Omega}+\Rho\tilde{\theta})\nonumber\\
& + k_{R}\tilde{e}_{\Omega} \cdot \tilde{e}_{R} + cC(R^{T}\hat{R})\tilde{e}_{\Omega} \cdot \tilde{e}_{\Omega}\nonumber\\
& + c\tilde{e}_{R} \cdot J^{-1}(-k_{R}\tilde{e}_{R}-k_{\Omega}\tilde{e}_{\Omega}+\Rho\tilde{\theta})\nonumber \\
& + \tilde{\theta}^{T}\Gamma^{-1}\dot{\tilde{\theta}}.\label{eq:Vdot2}
\end{align}

The resulting $\dot{V}$ can be separated into two parts, one with terms containing $\tilde{\theta}$ and the other without $\tilde{\theta}$,

\begin{align}
\dot{V} =& \dot{V}_{\dvone} + \dot{V}_{\dvtwo},\label{eq:V1dotpart}\\
\dot{V}_{\dvone} \coloneqq & \tilde{e}_{\Omega} \cdot(-\tilde{k}_{R}\tilde{e}_{R}-\tilde{k}_{\Omega}\tilde{e}_{\Omega})\nonumber \\
&+ \tilde{k}_{R}\tilde{e}_{\Omega} \cdot \tilde{e}_{R} + cC(R^{T}\hat{R})\tilde{e}_{\Omega} \cdot \tilde{e}_{\Omega} \nonumber \\
&+ c\tilde{e}_{R} \cdot J^{-1}(-\tilde{k}_{R}\tilde{e}_{R}-\tilde{k}_{\Omega}\tilde{e}_{\Omega}), \label{eq:V1dot}
\end{align}
\begin{align}
\dot{V}_{\dvtwo} \coloneqq & \tilde{e}_{\Omega} \cdot \Rho\tilde{\theta} + c\tilde{e}_{R} \cdot J^{-1}\Rho\tilde{\theta} + \tilde{\theta}^{T}\Gamma^{-1}\dot{\tilde{\theta}},  \label{eq:V2dot} \end{align}
%
Simplify \eqref{eq:V1dot} to get,
\begin{equation}
\begin{split}
\dot{V}_{\dvone} =& -k_{\Omega}\|\tilde{e}_{\Omega}\|^{2}- ck_{R}\tilde{e}_{R} \cdot J^{-1}\tilde{e}_{R}\\
&+ cC(\hat{R}^{T}R)\tilde{e}_{\Omega} \cdot \tilde{e}_{\Omega} - ck_{\Omega}\tilde{e}_{R} \cdot J^{-1}\tilde{e}_{\Omega}.
\end{split}
\end{equation}
This has the same form as \cite[(58)]{lee2011geo} and thus we have,
\begin{equation}
\dot{V}_{\dvone} \leq -\tilde{\eta}^{T}W\tilde{\eta},
\end{equation}
where 
\begin{align}
\label{eq:eta_def}
\tilde{\eta} = \begin{bmatrix}
\|{\tilde{e}_{R}}\|& \|{\tilde{e}_{\Omega}}\|
\end{bmatrix}^{T}
\end{align}
 and $W$ is,
\begin{equation}
W = \begin{bmatrix} \frac{ck_{R}}{\lambda_{M}(J)} & -\frac{ck_{\Omega}}{2\lambda_{m}(J)} \\
-\frac{ck_{\Omega}}{2\lambda_{m}(J)} & k_{\Omega}-c  \end{bmatrix},
\end{equation}
where $\lambda_{m}(A)$ and $\lambda_{M}(A)$ are the minimum and maximum eigenvalues of the matrix A, respectively. The constant $c$ is chosen such that $W$ is positive definite resulting in $\dot{V}_{\dvone} \leq 0$.

From \eqref{eq:V}, let $V_1$ be,
\begin{align}
V_1 = \frac{1}{2}\tilde{e}_{\Omega} \cdot J\tilde{e}_{\Omega}+\tilde{k}_{R}\tilde{\Psi}(\hat{R},R)+c\tilde{e}_{R} \cdot \tilde{e}_{\Omega}
\end{align}
such that $ V = V_1 + \frac{1}{2} \tilde \theta^T \Gamma^{-1} \tilde \theta$.  Then, as shown in detail in [13], we can show $V_1$ satisfies,
\begin{align}
\tilde{\eta}^T W_1 \tilde{\eta} \leq V_1 \leq \tilde{\eta}^T W_1 \tilde{\eta},
\end{align}
where, 
\begin{align}
W_1=\frac{1}{2}\begin{bmatrix}
\tilde{k}_R & -c \\ -c & \lambda_{m}(J)
\end{bmatrix},\,
W_2=\frac{1}{2}\begin{bmatrix}
\frac{2\tilde{k}_R}{2-\psi} & c\\ c & \lambda_{M}(J)
\end{bmatrix}
\end{align}
with,
\begin{align*}
\tilde{\Psi}(\hat{R}(t), R(t)) \leq \psi < 2,\quad \text{for any t},
\end{align*}
and $\tilde{\eta}$ defined in \eqref{eq:eta_def}. 

$V_1$, and $\dot{V}_{\dvone}$ are bounded as
\begin{gather}
\lambda_{m}(W_1)\| \tilde{\eta} \|^2 \leq V_1 \leq \lambda_{M}(W_2)\| \tilde{\eta} \|^2,\\
\dot{V}_{\dvone} \leq -\lambda_{m}(W)\|\tilde{\eta} \|^2,
\end{gather}
we can further show that,
\begin{gather}
\dot{V}_{\dvone}\leq -\beta V_1, \quad \beta = \frac{\lambda_{m}(W)}{\lambda_{M}(W_2)}. \label{eq:dV1_final}
\end{gather}
These results are used later in the proof.  Next, from \eqref{eq:V2dot}, we have,
\begin{align}
\dot{V}_{\dvtwo}&=\T{\theta}^T\big(\Rho^T\T{e}_{\Om} + c\Rho^TJ^{-T}\T{e}_R\big)  + \T{\theta}^T\Gamma^{-1}\dot{\T{\theta}}. \label{eq:V2dot_1}
\end{align}

Taking derivative of \eqref{eq:thetatilde} we have,
\begin{gather}
\dot{\T{\theta}}=\dot{\hat{\theta}}+\hatmap{\T{e}}_{\Om}(\hat{R}^TR)\theta - \hat{R}^TR\dot{\theta}.
\end{gather}

Therefore \eqref{eq:V2dot_1} can be updated as,
\begin{align}
\dot{V}_{\dvtwo}&=\T{\theta}^T\big[\Rho^T\T{e}_{\Om} + c\Rho^TJ^{-T}\T{e}_R\big] \nonumber \\
& \quad + \T{\theta}^T\Gamma^{-1}\big[\dot{\hat{\theta}}+\hatmap{\T{e}}_{\Om}(\hat{R}^TR)\theta - \hat{R}^TR\dot{\theta}\big]. \label{eq:V2dot_2}
\end{align}

 This is rewritten as,
\begin{align}
\dot{V}_{\dvtwo}&= \underbrace{\T{\theta}^T\big[\Rho^T\T{e}_{\Om} + c\Rho^TJ^{-T}\T{e}_R  + \Gamma^{-1}{\dot{\hat{\theta}}}\big]}_{\triangleq \dot{V}_{\dvtwo a}} \nonumber \\ 
&\quad + \underbrace{\T{\theta}^T\Gamma^{-1}\big[\hatmap{\T{e}}_{\Om}(\hat{R}^TR)\theta -\hat{R}^TR\dot{\theta}\big]}_{\triangleq \dot{V}_{\dvtwo b}}. \label{eq:V2dot_2parted}
\end{align}

We use a property of the $\Gamma-$projection operator as shown in \cite{lavretsky2011projection},
\begin{equation}
\label{eq:Proj}
\tilde{\theta}^{T}(\Gamma^{-1}Proj_{\Gamma}(\hat{\theta},y)-y) \leq 0,
\end{equation}
with the projection operator defined in \eqref{eq:adapt} and the projection function $y$ defined in \eqref{eq:y}, we have,
 
\begin{align}
\dot{V}_{\dvtwo a}\leq 0. \label{eq:dV2a}
\end{align} 
We assume that the uncertainty $\theta$ and its time derivative $\dot{\theta}$ are bounded. Furthermore, the projection operator in \eqref{eq:adapt} will also keep $\tilde{\theta}$ bounded (see \cite{Cao2006adapt} for a detailed proof about these properties.) Thus, we consider the following bounds,
\begin{align}
\norm{\T{\theta}} \leq \T{\theta}_{b} \quad \& \quad \norm{\dot{\theta}}\leq \dot{\theta}_b \quad \& \quad \norm{\theta} \leq \theta_b.
\end{align}

From \eqref{eq:V2dot_2parted}, $\dot{V}_{\dvtwo b}$ can be bounded as follows,
\begin{align}
\dot{V}_{\dvtwo b} &= \T{\theta}^T\Gamma^{-1}\big(\hatmap{\T{e}}_{\Om}(\hat{R}^TR)\theta -\hat{R}^TR\dot{\theta}\big)\nonumber \\
& \leq \norm{\T{\theta}}\norm{\Gamma^{-1}}\bigg( \norm{\T{e}_{\Om}}\norm{(\hat{R}^TR)}\norm{\theta} + \norm{(\hat{R}^TR)}\norm{\dot{\theta}}\bigg) \nonumber  \\
& \leq \T{\theta}_b\norm{\Gamma^{-1}}\big( \norm{\T{e}_{\Om}}{\theta}_b + {\dot{\theta}}_b\big).
\end{align}

From \eqref{eq:final_deOmtilde}, we can show that the ${\tilde{e}}_{\Omega}$ is decreasing for a right choice of $\tilde{k}_R$, $\tilde{k}_{\Omega}$ (since $\tilde{\theta}$ is bounded). Also, from \eqref{eq:eOmtilde_init}, $\tilde{e}_{\Omega}(0)$ is bounded. Initial value of $\tilde{e}_{\Omega}$ being bounded and $\tilde{e}_{\Omega}$ decreasing, implies $\tilde{e}_{\Omega}$ is bounded for all time. Therefore, let $\norm{\tilde{e}_{\Omega}}$ be bounded by $\tilde{e}_{\Omega b}$.  Then,
\begin{equation}
\dot{V}_{\dvtwo b} \leq \T{\theta}_b(\T{e}_{\Omega b}\theta_b + \dot{\theta}_b)\norm{\Gamma^{-1}},
\end{equation}

Choosing a large adaptation gain, $\Gamma$, would result in a very small $\Gamma^{-1}$ and thus the right side can be bounded to a small neighborhood $\delta$. Then,
\begin{align}
\dot{V}_{\dvtwo b} \leq \delta, \label{eq:dV2b}
\end{align}
and thus from \eqref{eq:dV2a} and \eqref{eq:dV2b}, we have,
\begin{align}
\dot{V}_{\dvtwo}\leq \delta. \label{eq:dV2_final}
\end{align}

Substituting \eqref{eq:dV1_final} and \eqref{eq:dV2_final} in \eqref{eq:V1dotpart} we get,
\begin{align}
\dot{V} & \leq -\beta V_1 + \delta \\
 & \leq -\beta( \underbrace{V_1 + \frac{1}{2}\tilde{\theta}^{T}\Gamma^{-1}\tilde{\theta}}_{V}) + \beta(\frac{1}{2}\tilde{\theta}^{T}\Gamma^{-1}\tilde{\theta}) + \delta \\
& \leq -\beta V + \beta\frac{{\T{\theta}}_b^2}{2}\norm{\Gamma^{-1}} + \delta
\end{align}
We finally have, 
\begin{equation}
\label{eq:Final_V_eq}
\dot{V} + \beta V \leq \beta \delta_V 
\end{equation}
with $\delta_V \triangleq \frac{{\T{\theta}}_b^2}{2}\norm{\Gamma^{-1}} + \frac{\delta}{\beta}$. In \eqref{eq:Final_V_eq}, if $V\geq \delta_V$ it results in $\dot{V}\leq 0$. As a result, by choosing a sufficiently large adaptation gain $\Gamma$, the Control Lyapunov Candidate function V, decreases exponentially to result in $V\leq \delta_V$, an arbitrarily small neighborhood $\delta_V$. 
As shown in \eqref{eq:Final_V_eq}, $V$ is an \textit{exponential Input-to-State Stable} Lyapunov function  \cite{kolathaya2018input}, and thus the attitude errors $(\tilde{e}_R,\tilde{e}_{\Omega})$ are exponential input-to-state stable.

\balance
	
	\bibliographystyle{IEEEtranS}

	\bibliography{references}

\begin{thebibliography}{10}
\providecommand{\url}[1]{#1}
\csname url@samestyle\endcsname
\providecommand{\newblock}{\relax}
\providecommand{\bibinfo}[2]{#2}
\providecommand{\BIBentrySTDinterwordspacing}{\spaceskip=0pt\relax}
\providecommand{\BIBentryALTinterwordstretchfactor}{4}
\providecommand{\BIBentryALTinterwordspacing}{\spaceskip=\fontdimen2\font plus
\BIBentryALTinterwordstretchfactor\fontdimen3\font minus
  \fontdimen4\font\relax}
\providecommand{\BIBforeignlanguage}[2]{{%
\expandafter\ifx\csname l@#1\endcsname\relax
\typeout{** WARNING: IEEEtranS.bst: No hyphenation pattern has been}%
\typeout{** loaded for the language `#1'. Using the pattern for}%
\typeout{** the default language instead.}%
\else
\language=\csname l@#1\endcsname
\fi
#2}}
\providecommand{\BIBdecl}{\relax}
\BIBdecl

\bibitem{ghaffar2016mrac}
A.~F. Abdul~Ghaffar and T.~S. Richardson, ``Position tracking of an
  underactuated quadrotor using model reference adaptive control,'' in
  \emph{AIAA Guidance, Navigation, and Control Conference}, 2016, p. 1388.

\bibitem{ackerman2016l1}
K.~Ackerman, E.~Xargay, R.~Choe, N.~Hovakimyan, M.~C. Cotting, R.~B. Jeffrey,
  M.~P. Blackstun, T.~P. Fulkerson, T.~R. Lau, and S.~S. Stephens, ``L1
  stability augmentation system for calspan's variable-stability learjet,'' in
  \emph{AIAA Guidance, Navigation, and Control Conference}, 2016, p. 0631.

\bibitem{cao2006design}
C.~Cao and N.~Hovakimyan, ``Design and analysis of a novel l1 adaptive
  controller, part i: Control signal and asymptotic stability,'' in
  \emph{American Control Conference}, 2006, pp. 3397--3402.

\bibitem{Cao2006adapt}
------, ``Design and analysis of a novel l1 adaptive controller, part ii:
  Guaranteed transient performance,'' in \emph{American Control Conference},
  2006, pp. 3403--3408.

\bibitem{monte2013adapt}
P.~De~Monte and B.~Lohmann, ``Position trajectory tracking of a quadrotor based
  on l1 adaptive control,'' \emph{at--Automatisierungstechnik}, vol.~62, no.~3,
  pp. 188--202, 2014.

\bibitem{dodenhoft2017design}
J.~Dodenh{\"o}ft, R.~Choe, K.~Ackerman, F.~Holzapfel, and N.~Hovakimyan,
  ``Design and evaluation of an l1 adaptive controller for nasa’s transport
  class model,'' in \emph{AIAA Guidance, Navigation, and Control Conference},
  2017, p. 1250.

\bibitem{goodarzi2013geometric}
F.~Goodarzi, D.~Lee, and T.~Lee, ``Geometric nonlinear pid control of a
  quadrotor uav on se (3),'' in \emph{2013 European Control Conference
  (ECC)}.\hskip 1em plus 0.5em minus 0.4em\relax IEEE, 2013, pp. 3845--3850.

\bibitem{hovakimyan2010l1}
N.~Hovakimyan and C.~Cao, \emph{ℒ1 Adaptive Control Theory: Guaranteed
  Robustness with Fast Adaptation}.\hskip 1em plus 0.5em minus 0.4em\relax
  SIAM, 2010.

\bibitem{kolathaya2018input}
S.~Kolathaya and A.~D. Ames, ``Input-to-state safety with control barrier
  functions,'' \emph{IEEE control systems letters}, vol.~3, no.~1, pp.
  108--113, 2018.

\bibitem{kulumani2016adapt}
S.~Kulumani, C.~Poole, and T.~Lee, ``Geometric adaptive control of attitude
  dynamics on so(3) with state inequality constraints,'' in \emph{American
  Control Conference}, 2016, pp. 4936--4941.

\bibitem{lavretsky2011projection}
E.~Lavretsky and T.~E. Gibson, ``Projection operator in adaptive systems,''
  \emph{arXiv preprint arXiv:1112.4232}, 2011.

\bibitem{lee2015adapt}
T.~Lee, ``Geometric adaptive control for aerial transportation of a rigid
  body,'' in \emph{Proceedings of the IMA Conference on Mathematics of
  Robotics}, 2015.

\bibitem{lee2011geo}
T.~Lee, M.~Leok, and N.~H. McClamroch, ``Geometric tracking control of a
  quadrotor uav on se(3),'' in \emph{Intl. Conference on Decision and Control},
  2010, pp. 5420--5425.

\bibitem{mahony2012multirotor}
R.~Mahony, V.~Kumar, and P.~Corke, ``Multirotor aerial vehicles,'' \emph{IEEE
  Robotics and Automation magazine}, vol.~20, no.~32, 2012.

\bibitem{mallikarjunan2012l1}
S.~Mallikarjunan, B.~Nesbitt, E.~Kharisov, E.~Xargay, N.~Hovakimyan, and
  C.~Cao, ``L1 adaptive controller for attitude control of multirotors,'' in
  \emph{AIAA Guidance, Navigation, and Control Conference}, 2012, p. 4831.

\bibitem{mellinger2011minimum}
D.~Mellinger and V.~Kumar, ``Minimum snap trajectory generation and control for
  quadrotors,'' in \emph{International Conference on Robotics and Automation},
  2011, pp. 2520--2525.

\bibitem{michini2009adapt}
B.~Michini, ``Modeling and adaptive control of indoor unmanned aerial
  vehicles,'' Ph.D. dissertation, Massachusetts Institute of Technology, 2009.

\bibitem{quan2015bi}
Q.~Nguyen and K.~Sreenath, ``L1 adaptive control for bipedal robots with
  control lyapunov function based quadratic programs,'' in \emph{American
  Control Conference}, 2015, pp. 862--867.

\bibitem{quan2017introduction}
Q.~Quan, \emph{Introduction to multicopter design and control}.\hskip 1em plus
  0.5em minus 0.4em\relax Springer, 2017.

\bibitem{koushil2013geo}
K.~Sreenath, T.~Lee, and V.~Kumar, ``Geometric control and differential
  flatness of a quadrotor uav with a cable-suspended load,'' in \emph{Intl.
  Conference on Decision and Control}, 2013, pp. 2269--2274.

\bibitem{whitehead2010mrac}
B.~Whitehead and S.~Bieniawski, ``Model reference adaptive control of a
  quadrotor uav,'' in \emph{AIAA Guidance, Navigation, and Control Conference},
  2010, p. 8148.

\bibitem{wu2014geo}
G.~Wu and K.~Sreenath, ``Geometric control of multiple quadrotors transporting
  a rigid-body load,'' in \emph{Intl. Conference on Decision and Control},
  2014, pp. 6141--6148.

\bibitem{zuo2014adapt}
Z.~Zuo and P.~Ru, ``Augmented l1 adaptive tracking control of quad-rotor
  unmanned aircrafts,'' \emph{IEEE Transactions on Aerospace and Electronic
  Systems}, vol.~50, no.~4, pp. 3090--3101, 2014.

\end{thebibliography}
	
\end{document}